\documentclass{elsarticle}

\journal{Journal of the Franklin Institute}

\usepackage{caption}

\usepackage{subcaption}

\usepackage{amsfonts,amsmath,amssymb}
\usepackage{graphicx}

\renewcommand{\geq}{\geqslant}
\renewcommand{\leq}{\leqslant}

\usepackage{siunitx}

\newtheorem{remark}{Remark}

\newtheorem{lemma}{Lemma}
\newtheorem{theorem}{Theorem}

\usepackage{color}

\begin{document}

\begin{frontmatter}

\title{On a nodal observer for 
 a semilinear model for the flow in 
 gas networks\tnoteref{t1}}
\tnotetext[t1]{This work was
supported by DFG
in the framework of the Collaborative Research Centre
CRC/Transregio 154,
Mathematical Modelling, Simulation and Optimization Using the Example of Gas Networks,
Project C03 and  Project  C05.}

\author{Martin Gugat\corref{cor1}}

\address{Friedrich-Alexander-Universit\"at Erlangen-N\"urnberg (FAU),
Department of Data Science, Cauerstr. 11, 91058
Erlangen, Germany}

\ead{martin.gugat@fau.de}

\author{Jan Giesselmann}

\author{Teresa Kunkel}
\address{Technische Universit\"at Darmstadt,
Fachbereich Mathematik, Dolivostr. 15, 64293 Darmstadt, Germany.}

\cortext[cor1]{Corresponding author}



\begin{abstract}
The flow of gas through networks of pipes
can be modeled by coupling 
hyperbolic systems of 
partial differential equations
that describe the flow through
the pipes that form the edges of
the graph of the network 
by algebraic node conditions
that model the
flow through the vertices of
the  graph.
In the network,
measurements of the state 
are available at certain points in space.
Based upon these nodal observations,
the complete system state can be 
approximated using an observer system. 
In this paper we present
a nodal observer,
and prove that the
state of the observer system
converges to the original state
exponentially fast.
 Numerical experiments confirm the theoretical findings. 
\end{abstract}

\begin{keyword}
Network, node conditions, gas transportation network,
observability inequality,
exponential synchronization,
networked hyperbolic system,
semilinear hyperbolic pde
\MSC[2010] 35L04 \sep 49K20
\end{keyword}

\end{frontmatter}





This work was
supported by DFG
in the framework of the Collaborative Research Centre
CRC/Transregio 154,
Mathematical Modelling, Simulation and Optimization Using the Example of Gas Networks,
Project C03 and  Project  C05.

\section*{Introduction}
In this contribution we study the
problem to construct an observer for   the flow of gas through
networks of pipelines.
In this application,
in general, many  pipelines are very long
and the graph of the network can
be quite complex.
We consider a
semilinear model for a gas pipeline network where at the nodes the solutions
for the adjacent pipes are coupled by  algebraic node conditions that require
the conservation of mass and the continuity of the pressure.
The semilinear model is 
a resonable simplification of the
quasilinear isothermal Euler equations  
if  the Mach number of the flow
is  quite small.
In the operation of gas transportation systems this is the case, since 
the velocity of the gas flow is much smaller than the
sound speed.
The eigenvalues of the quasilinear system have the form
$c + v $ and $-c + v$, where $v$ denotes
the velocity of the gas and $c$ denotes the sound speed.
In order to obtain a semilinear model
the eigenvalues are replaced by the sound speed,
that is by $-c$ and $c$, at some reference density.
%
While we linearize the convection,
in our semilinear model 
we keep the nonlinear  source term.
This source term plays an essential role in the
model of gas network flow, since
the friction effects lead to a
decrease of the pressure along each pipe in
the direction of the flow.

In this paper, 
we present a nodal  observer
system for the network flow
where the coupling to
the original system 
is governed 
at each vertex $v$
by a parameter $\mu^v\in [-1,1]$.
We show that the observer system yields
an approximation of the
state in the original system 
where the error decays 
exponentially fast in
the sense of the $L^2$-norm 
if the $L^\infty$-norm of the initial state is sufficiently small.
Moreover, we also show that
the $H^1$-norm of the approximation error decays exponentially fast if the
  state of
  the original system is sufficiently 
  regular (i.e. in $W^{1,\infty}$).
It is
desirable to have such an observer since
it allows to obtain a reliable 
approximation for the complete
state in the original system.
%
The state estimation that is obtained
with the observer system can be used in the construction of feedback laws.
The proofs are
based upon
nodal observability inequalities
for the flow.

%
%

Observers, using distributed measurements, have been constructed for semi-linear hyperbolic equations \cite{CCDB12,CIM15} and quasi-linear hyperbolic equations \cite{BMP15}.
The results for semi-linear hyperbolic equations are based on stabilization results for (locally) damped wave equations \cite{HST12}. The result for quasi-linear problems is based on a kinetic formulation and, thus, uses a rather different approach.
In \cite{HAK16} the backstepping method has been used to construct 
a boundary observer for semi-linear hyperbolic problems.
For results about the recovery
of an unknown initial state
using an observer see \cite{IMT20}.
In \cite{FERRANTE}, the design of
boundary observers  for a linear system of ODEs in cascade with hyperbolic PDEs is studied
and more references on  observer design are given.
We want to emphasize that 
the novelty of our contribution 
is the construction of 
an observer
for a  system that is governed by 
 networked semilinear pdes and uses 
observations that are located pointwise in  space, whereas in the previous contributions distributed observations 
coming  from subdomains
in space have been considered.




%

This paper has the following structure.
In Section \ref{isothermalEuler} we introduce
the quasilinear isothermal Euler equations.
In Section  \ref{sec:Riemann}, we present the
corresponding Riemann invariants and  transform the system
in diagonal form.
In Section \ref{sec:networks} we present
the node conditions that model the flow through
the junctions in  a gas  pipeline network.
 Then we derive the semilinear model
that provides an approximation for small gas velocities.

We are working in the framework of
solutions that are defined
through a fixed point iteration
along the characteristics.
The definition of the fixed point mapping
is derived from the integral
equations along the characteristic curves
that are known a priori for the semilinear model. A well--posedness result that is
based on this approach is presented in 
Section \ref{wellposedness}.

In Section \ref{Sec:Exponential_decay} 
we show that the error of the observer 
converges  to zero
in $L^2$
exponentially fast.
For the proof,
in Section \ref{l2observability} we first introduce a quadratic
$L^2$-Lyapunov function 
and show  that it decays exponentially 
fast on finite time intervals without
additional constraints on the lengths of the pipes.
In the proof we use an observability inequality for the $L^2$-norm.
Then in Section \ref{observabilityh1}
we define a  quadratic Lyapunov function
with exponential weights
to show that the time derivatives also  decay exponentially fast.
This yields the exponential decay of the $H^1$-norm
of the state for initial states with sufficiently small $H^1$-norm.
This is also shown in 
Section \ref{Sec:Exponential_decay}.
However, the result about the exponential
decay in $H^1$ requires strict
bounds on the  lengths of the pipes.
In Section \ref{Sec:Numerical}
at the end of the paper, numerical
results are presented.


\section{Motivating the semilinear model from the 
Euler equations}
\label{motivation}
\subsection{The
isothermal Euler equations along  the pipes}
\label{isothermalEuler}
\noindent
The isothermal Euler equations
as a model for the flow through gas pipelines
have  been stated for example  in  \cite{bandahertyklar}.
Let a  finite  graph $G=(V,\,E)$ of a pipeline network be given.
Here $V$ denotes the set of vertices and $E$ denotes the set of edges.
Each edge $e\in E$  corresponds to  an interval $[0,\, L^e]$ that represents
a pipe of length $L^e>0$.
Let $D^e>0$ denote the diameter  and
{
$\lambda_{fric}^e
>
0$
}
the  friction coefficient in pipe $e$.
Define
$\theta^e = \frac{\lambda_{fric}^e}{D^e}$.
%
Let $\rho^e$ denote the gas density, $p^e$ the pressure and $q^e$ the mass flow rate.
Since we will linearize the convective part of the equations, we will not specify any pressure law. The isothermal or isentropic Euler equations are hyperbolic provided the pressure $p=p(\rho)$ is given as a monotone increasing function of the density. We will assume that this monotonicity is strict.
Typical examples are isentropic law $p(\rho)= a\rho^\gamma$ with $a >0$ and $\gamma > 1$ and the model of the American Gas Association (AGA), see \cite{gugatulbrich},
\[ p(\rho)= \frac{ R_sT \rho}{1- \tilde \alpha \rho}\]
where $T$ is the temperature of the pipe, $R_s$  is the gas constant and $\tilde \alpha \leq 0$.
Note that for $\tilde \alpha=0$ the AGA model reduces to the isothermal law $p(\rho)=R_sT \rho$.
We study a model that is based upon the $2 \time 2$ Euler equations, i.e. there is no energy conservation equation,
\begin{equation}
\label{isothermaleuler}
\left\{
\begin{array}{l}
\rho_t^e + q_x^e = 0,
\\
q_t^e + \left( p^e + \frac{(q^e)^2}{\rho^e} \right)_x = - \frac{1}{2} \theta^e \frac{q^e \; |q^e|}{\rho^e}
\end{array}
\right.
\end{equation}
that govern the flow through  pipe $e$.
In  order to deal with large networks,
it is desirable to 
replace the quasilinear Euler equations
by a simpler semilinear model.
In the  sequel, we will present such a model.
We proceed in the following way.
First we transform the system
to a system in diagonal form
for the Riemann invariants.
Then we observe that in the
normal operational of range gas pipeline
flow (that is for slow transients), the matrix in the diagonal
form is close to a constant matrix.
We 
obtain a semilinear model by
replacing the system matrix by 
this constant matrix.



\subsection{The system in terms of Riemann invariants}
\label{sec:Riemann}
Hyperbolic systems of conservation laws can be diagonalized  by Riemann invariants, i.e. quantities whose derivatives with respect to the conserved variables are orthogonal to all but one of the left eigenvectors of the Jacobian of the flux.

Here  the flux, its Jacobian, and the eigenvectors are given by
\[ f(\rho,q)= \begin{pmatrix} q \\ \frac{q^2 }{\rho} + p\end{pmatrix}, \
Df(\rho,q)= \begin{pmatrix} 0 & 1 \\ - \frac{q^2 }{\rho^2} + p'(\rho) & 2 \frac{q }{\rho} \end{pmatrix},\
\ell_\pm(\rho,q) = \begin{pmatrix}  -1 \\ \frac{q}{\rho} \pm \sqrt{p'(\rho)} .
\end{pmatrix}
\]
As explained in \cite[Chapter 7.3]{dafermosbook} every $2 \times 2$ system of hyperbolic conservation laws is endowed with a system of Riemann invariants. For the Euler equations with general pressure law they are given by
\[ R_\pm (\rho,q)=  \tilde R(\rho) \pm \frac{q}{\rho}\]
where $\tilde R$ is defined  by
$
\tilde R(\rho)=    \int_1^\rho \frac{\sqrt{p'(r)}}{r}   \, dr .
$
Specific formulas for Riemann invariants  for the isentropic law $p(\rho)=a \rho^\gamma$ can be found in \cite[Chapter 7.3]{dafermosbook}
while the specific formulas for the AGA model was computed in 
 \cite{gugatulbrich}.
 Note that Riemann invariants are only unique up operations that leave the direction of the gradient unchanged.

\section{The Node Conditions for the Network Flow}
\label{sec:networks}
In this section we introduce the  coupling conditions
that model the flow through the nodes of the network.
%
%
For any node $v \in V$ let $E_0(v)$ denote the set of edges in the graph that are incident to
$v\in V$
and let $x^e(v)\in \{0,L^e\}$ denote the
end of the interval~$[0,L^e]$ that corresponds
to the edge $e$ that is adjacent to $v$.
Let $V_0(e)$ denote the set of nodes adjacent to some edge $e$.
Define
\begin{equation}
\label{mathfraksdefinition}
{\mathfrak s}(v, \, e):=\left\{
\begin{array}{rll}
-1 & {\rm if} & x^e(v)=0 \;\mbox{\rm and }\; e\in E_0(v),
\\
1 & {\rm if} & x^e(v)=L^e \;\mbox{\rm and }\; e\in E_0(v),
\\
0 &{\rm if}  &  e\not\in E_0(v).
\end{array}
\right.
\end{equation}

We impose the
Kirchhoff condition
\begin{equation}
\label{Kirchhoff}
\sum_{e \in E_0(v)}
{\mathfrak s}(v, \, e)
\,
 (D^e)^2\,
q^e(x^e(v))
=  
0
\end{equation}
that expresses conservation of mass at the nodes.

In order to close the system, additional coupling conditions are needed .
A typical choice, leading to well-posed Riemann problems \cite{bandahertyklar}, is 
to require 
the  continuity of the pressure at $v$,
which   means that
 for all $e$, $f\in E_0(v)$   we  have  
 \begin{equation}
\label{pressurecontinuity}
p(\rho^e(t,x^e(v))) = p(\rho^f(t,x^f(v))).
\end{equation}
another choice, which was advocated by \cite{reigstad} is continuity of enthalpy:
 for all $e$, $f\in E_0(v)$   we  have  
 \begin{equation}
\label{enthalpycontinuity}
P'(\rho^e(t,x^e(v)))  + \frac{(q^e(t,x^e(v))^2}{(\rho^e(t,x^e(v))^2} = P'(\rho^e(t,x^f(v)))
+ \frac{(q^f(t,x^e(v))^2}{(\rho^f(t,x^e(v))^2}
\end{equation}
where $P=P(\rho)$ is the pressure potential that is defined  by
\begin{align} \label{eq:sys3a}
P(\rho) = \rho \int_1^\rho \frac{p(r)}{r^2} dr,
\end{align}

Our interest is in simplified models where velocities are much smaller than the speed of sound. Let us note that in the $\tfrac{q}{\rho} \rightarrow 0$ limit both \eqref{pressurecontinuity} and \eqref{enthalpycontinuity} enforce continuity of densities, since $p(\rho)$ and $P'(\rho)$ are both strictly monotone increasing.  Since both conditions coincide (asymptotically) in the limit of interest, we will use \eqref{pressurecontinuity} for our further discussions, irrespective of the question of which coupling is correct for general flows.

Now also state the node conditions in terms of Riemann invariants.

Define the vectors 
$R_{in }^v(t)$,
$R_{out}^v(t)
\in {\mathbb R}^{|E_0(v)|}
$
in the following manner:
\\
If $e  \in  E_0(v)$ and
${\mathfrak s}(v,\, e) =1$,  $R^e_+(t, \, x^e(v))$
is a component of
 $R_{in }^v(t)$
 that we refer to as
 $R_{in }^e(t,x^e(v))$
and if
 $e  \in  E_0(v)$ and
${\mathfrak s}(v,\, e) =-1$, $R_{in }^v(t)$ contains $R^e_-(t,\, x^e(v))$
as a component
that we also refer to as
 $R_{in }^e(t,x^e(v))$.
\\
Moreover, if $e  \in  E_0(v)$ and
${\mathfrak s}(v,\, e) =1$,  $R^e_-(t,\,x^e(v))$
is a component of
 $R_{out}^v$
 that we refer to as
 $R_{out }^e(t,x^e(v))$
and if
 $e  \in  E_0(v)$ and
${\mathfrak s}(v,\, e) =-1$, $R_{out}^v(t)$ contains $R^e_+(t, x^e(v))$
as a component
that we also refer to as
 $R_{ out }^e(t,x^e(v))$.

We assume that the components
are ordered in such a way that
the $j$-th component of
$R_{out}^v$
corresponds to the same edge
$e\in E_0(v)$
as the 
$j$-th component of
$R_{in}^v$.


\begin{lemma}For any node $v\in V$ of the  graph 
and $e\in E_0(v)$ the node conditions
(\ref{pressurecontinuity}),
(\ref{Kirchhoff})
can be written in the form of the  linear equation
\begin{equation}\label{Omegav}
R_{out}^e(t,\, x^e(v)) = - R^e_{in}(t,\, x^e(v) ) + \omega_v  \sum\limits_{g\in E_0(v)} (D^g)^2 \, R^g_{in}(t,\, x^g(v))
.
\end{equation}
where
\begin{equation}
\label{omegavdefinition}
 \omega_v :=  \frac{2}{\sum\limits_{ f \in E_0(v) } (D^f)^2 }.
\end{equation}
\end{lemma}
\textbf{Proof:} 
Equation (\ref{Omegav}) implies
that for all $e\in E$, the value of $R_+^e(t,\, x^e(v)) + R^e_-(t,\, x^e(v)) $ is
the same, which implies that the value of $\tilde R(\rho^e)$ is
independent of $e$. Since $\rho \mapsto \tilde R(\rho)$ is strictly monotone increasing, this implies  (\ref{pressurecontinuity}).
%
Moreover, (\ref{Omegav}) implies
\[
\sum_{e \in E_0(v)} (D^e)^2 \left[ R^e_{out}(t,\, x^e(v)) - R^e_{in}(t,\, x^e(v)) \right] = 0.
\]
Due to (\ref{pressurecontinuity}) this implies that  equation (\ref{Kirchhoff}) holds.

%

For each   $v\in  V$, let a number
$\mu^v\in [-1,\, 1]$ be given.

For a boundary node
$v\in V$ where $|E_0(v)|=1$ 
we state the boundary  conditions
in terms of Riemann invariants in the form
\begin{eqnarray}
\label{rbriemann2neu}
R^e_{out}(t,\, x^e(v)) & = 
( 1 - \mu^v) \, u^e(t) + \mu^v R^e_{in}(t,\, x^e(v)).
\end{eqnarray}

\begin{remark}
For $\mu^v=0$, \eqref{rbriemann2neu} encodes Dirichlet boundary conditions for the incoming Riemann invariant. For $\mu^v=1$ \eqref{rbriemann2neu}
expresses a Dirichlet boundary condition for the 
velocity,
and for
$\mu^v = -1$
we obtain a Dirichlet boundary condition for the pressure.

\end{remark}


For $e\in E$, let  $\nu^e   =  \frac{1}{4 }\theta^e $ be given and define
\begin{equation}
\label{sigmadefinition}
\sigma^e (R_+^e,\,R_-^e)
=
\nu^e  \, \left|R_+^e - R_-^e\right| \, (R_+^e - R_-^e).
\end{equation}
Define $\tilde \Delta^e$ as the  diagonal $2 \times 2$ matrix
that contains the eigenvalues
\begin{eqnarray}
\label{lambdapmgleichung}
\tilde \lambda^e_\pm
&
=
&
\frac{R_+^e - R_-^e}{2}  \pm \sqrt{p'\left(\tilde R^{-1}   \left(\frac{R_+^e + R_-^e}{2}\right) \right)  }
\end{eqnarray}
In terms of the Riemann invariants, the quasilinear  system (\ref{isothermaleuler}) has the following diagonal form:
\[
\partial_t 
\left(
\begin{array}{r}
R_+^e
\\
R_-^e
\end{array}
\right)
 +
\tilde  \Delta^e
\partial_x
\left(
\begin{array}{r}
R_+^e
\\
R_-^e
\end{array}
\right)
 =
\sigma^e (R_+^e,\,R_-^e)
\;
\left(
\begin{array}{r}
-1
\\
1 \end{array}
\right)
\]
With a reference density
$\rho_{\text{ref}}>0$,
define the number
$c = \sqrt{p'(\rho_{ref})}$.

In order to  simplify the model, we replace the  eigenvalues by
the constants 
\begin{eqnarray}
\label{lambdapmgleichungkonstant}
\lambda^e_\pm  &  = &
\pm c,
\end{eqnarray}
This definition implies that $\lambda^e_- = - \lambda^e_+$.
Moreover, for all $e,\, f \in E$  we have
$\lambda^e_+ = \lambda^f_+$.
Define $\Delta^e$ as the  diagonal $2 \times 2$ matrix
that contains the eigenvalues  $\lambda^e_+$ and   $\lambda^e_-$.
The approximation of
$\tilde \lambda^e_+$  by $\lambda^e_+$
and $\tilde \lambda^e_-$  by $\lambda^e_-$  is justified by
the fact that in the practical applications, the fluid velocity is several meters per second while the speed of sound is several hundred meters per second, i.e., $v$ can be neglected relative to $c$. In addition, the variation of the speed of sound due to density variations is rather small.  In contrast, the friction term cannot be neglected as this would cause a large \textit{relative} error.
In this way, we obtain a semilinear model. We do not claim that solutions to the isothermal Euler equations and the semilinear system are close to each other for all times, but we do expect that solutions to both systems share important qualitative features.
Let us note that the difference between the models becomes smaller the closer solutions get to equilibrium.
This is in agreement with the simplifications that we have made for the coupling conditions, as discussed below equation \eqref{enthalpycontinuity}.
%
%
With the diagonal  matrix $\Delta^e$, 
the  semilinear model has the following form:
\[
{\bf (S)}
\left\{
\begin{array}{l}
S_+^e(0, x)     =  y_+^e(x),\,   x\in (0,L^e),  \, e\in E,\\
S_-^e(0, x)     =  y_-^e(x), \,  x\in (0,L^e),  \, e\in E,\\
S_{out}^e(t, x^e(v))   = (1- \mu^v) \,  u^e(t)  + \mu^v S_{in}^e(t, L^e)  , \,  t\in (0,T),\;   \mbox{\rm if } \;|E_0(v)|=1,
\\ 
S_{out}^e(t, x^e(v))     =
 - S^e_{in}(t,\,x^e(v)) + \omega_v \sum\limits_{g\in E_0(v)} (D^g)^2 \, S^g_{in}(t,\, x^g(v))
 ,\,  t\in (0,T), 
 \\
\hfill  \;   \mbox{\rm if } \;|E_0(v)|\geq 2,
\\
 \partial_t\!
\left(
\begin{array}{r}
S_+^e
\\
S_-^e
\end{array}
\right)
 +
 \Delta^e
\partial_x\!
\left(
\begin{array}{r}
S_+^e
\\
S_-^e
\end{array}
\right)
 =
\sigma^e (S_+^e,\,S_-^e)
\;
\left(
\begin{array}{r}
-1
\\
1 \end{array}
\right)
\;\;
\mbox{\rm on }\; [0,T]\times [0,L^e],\, e\in E.
\end{array}
\right.
\]

Note that for any given $(S^e_+,\, S^e_-) \in \mathbb{R}^2$ we have
\[ p^e= p\left(\tilde R^{-1}\left( \frac{S^e_+ + S^e_-}{2} \right) \right) \]
For the special case of the isentropic and AGA model we have

\[p^e=  a^{-\frac{1}{2\gamma} }\left( \frac{S^e_+ + S^e_-}{2} + a^{\frac12}\right)^{\frac{2\gamma}{\gamma -1}}  \quad \text{ and } \quad p^e = \exp\left( \frac{S^e_+ + S^e_-}{2} \right),\]
respectively.
In particular, this implies $p^e>0$.
On account of the physical
interpretation of the pressure  it is very desirable that for
the solutions we have $p^e>0$.
This is an advantage of the model that is given by
system ${\bf (S)}$.

A similar  semilinear model for gas transport 
has been studied in \cite{hintermuellerstrogies} in the context of
identification problems.
The model in  \cite{hintermuellerstrogies} 
has the disadvantage that the matrix of the linearization of
the source term is indefinite. However, the results 
from  \cite{hintermuellerstrogies} 
can
be adapted to the model that we consider in this paper.

We introduce the  observer system  ${\bf (R)} $
that depends on numbers
$\mu^v\in  [-1,1]$
that are given for all $v\in V$
and control
the flow of information from
the original system to the observer system.
For an interior node with $\mu^v=0$,
the values 
at the node $v$ in the observer system
are fully determined by
the information from the original system.

The observer 
system  has a similar
structure as ${\bf (S)}$: 
\[
{\bf (R)}
\left\{
\begin{array}{l}
R_+^e(0, x)    =  z_+^e(x),\,   x\in (0,L^e),  \, e\in E,\\
R_-^e(0, x)     =  z_-^e(x), \,  x\in (0,L^e),  \, e\in E,\\
R_{out}^e(t, x^e(v))   =  (1-\mu^v) \, u^e( t )  + \mu^v  \, R_{in}^e(t, x^e(v))  
  , \,  t\in (0,T),\;   \mbox{\rm if } \;
  |E_0(v)|=1,
\\ 
R_{out}^e(t, x^e(v))     =
 %
  S^e_{out}(t,\, x^e(v))
 -
 \mu^v \left[ R_{in}^e(t, x^e(v))  -  S^e_{in}(t,\, x^e(v))  \right]
 \\
\hfill + \mu^v \,  \omega_v \sum\limits_{g\in E_0(v)} (D^g)^2 \, 
 \left[  R^g_{in}(t,\, x^g(v))  -   S^g_{in}(t,\, x^g(v))  \right]
 ,
 \\
\hfill t\in (0,T), 
 \;   \mbox{\rm if } \;
  |E_0(v)|\geq 2,
\\
\\
 \partial_t \!
\left(
\begin{array}{r}
R_+^e
\\
R_-^e
\end{array}
\right)
 +
 \Delta^e
\partial_x
\left(
\begin{array}{r}
R_+^e
\\
R_-^e
\end{array}
\right)
 =
\sigma^e (R_+^e,\,R_-^e)
\left(
\begin{array}{r}
-1
\\
1 \end{array}
\right)
\mbox{\rm on }\; [0,T]\times [0,L^e],\, e\in E.
\end{array}
\right.
\]
The initial state
$(z_+^e, z_-^e)$
represents an estimation of
the initial state of the original system.
Note that the data  $S^e_-(t,\, 0)$ from the original system
${\bf (S)}$ 
enter the system state
through the node conditions.
Indeed, for $\mu^v=0$ the Riemann invariants of the observer coincide with the Riemann invariants of the observed system. In contrast, for $\mu^v=0$ the observer satisfies the same coupling conditions as the observed system, i.e. no measurement information is inserted.

For the analysis of the exponential decay, we study the difference 
\[\delta = R - S\]
between the state $R$ that is
generated in the observer ${\bf (R)} $ and the original state $S$. 
For the difference  $\delta$ we obtain the system
\[
{\bf (D{iff})}
\left\{
\begin{array}{l}
\delta_+^e(0, x)     =  z_+^e (x) -  y_+^e(x) ,\,   x\in (0,L^e),  \, e\in E,\\
\delta_-^e(0, x)     =  z_-^e(x) - y_-^e(x)  , \,  x\in (0,L^e),  \, e\in E,\\
\delta_{out}^e(t, x^e(v))   =  \mu^v \, \delta_{in}^e(t, x^e(v)) 
  , \,  t\in (0,T),\;   \mbox{\rm if } \;  |E_0(v)|=1,
\\ 
\delta_{out}^e(t, x^e(v))     =
 - \mu^v \,  \delta^e_{in}(t,\, x^e(v))  
 + \mu^v \,  \omega_v \sum\limits_{g\in E_0(v)} (D^g)^2 \, 
    \delta^g_{in}(t,\, x^g(v) ),   
 \\
 \hfill \,  t\in (0,T), 
 \;   \mbox{\rm if } \;  |E_0(v)| \geq 2,
\\
 \partial_t
\left(
\begin{array}{r}
\delta_+^e
\\
\delta_-^e
\end{array}
\right)
 +
 \Delta^e
\partial_x
\left(
\begin{array}{r}
\delta_+^e
\\
\delta_-^e
\end{array}
\right)
 =
\left(
\sigma^e (\delta_+^e+S_+^e,\,\delta_-^e + S_-^e)
-
\sigma^e (S_+^e,\,  S_-^e)
\right)
\left(
\begin{array}{r}
-1
\\
1 \end{array}
\right)
\\
\;\;
\mbox{\rm on }\; [0,T]\times [0,L^e],\, e\in E.
\end{array}
\right.
\]
For the values of the solutions
of {\bf (D{iff})} at the nodes,
we have the following lemma.
\begin{lemma}
Let $\{\delta_\pm^e\}_{e \in E}$ be a solution of {\bf (D{iff})}. Then,
at any node $v\in V$   
\begin{equation}
    \label{middlenode_pl}
\sum_{e\in E_0(v)}
  (D^e)^2 \, |\delta^e_{out}(\tau,\, x^e(v)) |^2 
= |\mu^v|^2 \sum_{e\in E_0(v)}
  (D^e)^2 \,   |\delta^e_{in}(\tau,\, x^e(v))|^2, 
  \end{equation}
\end{lemma}

\textbf{Proof:}
 We only give the proof for inner nodes since the result is straightforward for boundary nodes, i.e., for $|E_0(v)|=1.$
 Let $v$ be arbitrary but fixed. 
Let us note that
\begin{multline}\label{eq:interm0b}
\sum_{e\in E_0(v)}
  (D^e)^2 \,\big[  \delta^e_{out}(\tau,\, x^e(v))
  - \mu^v   \delta^e_{in}(\tau,\, x^e(v)) \big]\\
  =\sum_{e\in E_0(v)}  (D^e)^2 \,\big[  
  - 2 \mu^v   \delta^e_{in}(\tau,\, x^e(v))   +  \mu^v \,  \omega_v \sum\limits_{g\in E_0(v)} (D^g)^2 \, 
    \delta^g_{in}(t,\, x^g(v) )\big]=0
\end{multline}
where  we have used the definition of $\omega_v$ in the last step.
We multiply \eqref{eq:interm0b} with $ \delta_{out}^e(t, x^e(v))     
 + \mu^v \,  \delta^e_{in}(t,\, x^e(v)) $ and obtain
\begin{equation}\label{eq:interm1}
\sum_{e\in E_0(v)}
  (D^e)^2 \,\big[  |\delta^e_{out}(\tau,\, x^e(v)) |^2 
  - |\mu^v|^2  |\delta^e_{in}(\tau,\, x^e(v)) |^2 \Big]=0
\end{equation}
where we have used that 
$ \delta_{out}^e(t, x^e(v))     
 + \mu^v \,  \delta^e_{in}(t,\, x^e(v)) $
 is independent of $e$. 
 Equation \eqref{eq:interm1} is equivalent to the statement of the lemma by elementary operations.



\section{A well--posedness result}
\label{wellposedness}
%
%
%
In the semilinear model that we consider, 
the constant eigenvalues in the diagonal system matrix
define two families of
characteristics with constant slopes $c$ and $-c$.
For $e\in E$, define the sets
\[
\Gamma_+^e 
=
\{0\} \times [0,\, L^e] \cup [0,\, T] \times \{0\},
\;\;\;
\Gamma_-^e 
=
\{0\} \times [0,\, L^e] \cup [0,\, T] \times \{L^e\}.
\]
For $t \geq 0$, $ e \in E$
and the space variable $x \in [0, \, L^e]$ 
we define  the ${\mathbb R}^2$-valued function
$\xi^e_\pm(s,\, x,\, t)$ as the solution of the initial value problem
\[
\left\{
\begin{array}{r}
\xi^e_\pm(t,\, x,\, t)= (t,\, x),
\\
\partial_s \xi^e_\pm(s,\, x,\, t) = ( 1,\, \pm c).
\end{array}
\right.
\]
This implies that
\[
\xi^e_+(s,\, x,\, t)= (s,\, x + c (s-t)),
\quad
\xi^e_-(s,\, x,\, t)= (s,\, x - c (s-t)).
\]
Define the points
\[
P_0^{e\pm}(t,\, x) = \Gamma^e_\pm\cap\{\xi^e_\pm(s,\, x,\, t), \, s\in \mathbb{R}\}\in {\mathbb R}^2
.
\]
For the $t$-component of $P_0^{e\pm}(t,\, x)$ we use the notation 
$t^e_{\pm}(x, \, t)\geq 0$.
For the discussion of
the well--posedness we focus
on the discussion of ${\bf (S) }$.
The discussion for the
observer system
${\bf (R) }$
and  the error system
${\bf (Diff) }$ is analogous.
The solution of  ${\bf (S) }$   can be defined by 
rewriting the partial differential equation in the system
as   integral equations along these characteristic curves, that is
\begin{equation}
\label{integralequation}
S^e_\pm(t,\, x)=
S^e_\pm(P_0^{e\pm}(t,\, x))
-\int_{t^e_{\pm}(x,\, t)}^t
 \pm\sigma^e(S^e_+,\, S^e_-)
(\xi^e_\pm(s,\, t,\, x))\, ds.
\end{equation}
Note that 
almost everywhere the values of
$S^e_\pm(P_0^{e\pm}(t,\, x))$
are given
on
$\Gamma^e_\pm$
either by the initial data,
that is $y^e_+$, $y^e_-$ respectively 
(if the $t$-component of $P_0^{e\pm}(t,\, x)$,
that is
$t^e_{\pm}(x, \, t)$  is zero),
the boundary condition
(\ref{rbriemann2neu}) 
(if the $x$-component of $P_0^{e\pm}(t,\, x)$ is zero or $L^e$
and for the corresponding node
we have $|E_0(v)|=1$)
or  by the node condition
(\ref{Omegav})  if
for the corresponding node
we have $|E_0(v)|\geq 2$.
For a finite time interval $[0, T]$,
the characteristic curves that start at $t=0$ 
with the information from the initial data 
reach a point at the terminal time after a finite number of reflections
at the boundaries  $x=L^e$ ($e\in E$) or  $x=0$.

The definition of the solutions
of semilinear hyperbolic boundary value problems
based upon  (\ref{integralequation}) is described for example in \cite{brokate}. 
For $L^\infty$-solutions, we have the following theorem.

\begin{theorem}
\label{existence}
Let $T>0$,
a real number $J$  and a 
number $M>0$ be given.

Then there exists a number $\varepsilon(T, \, M)>0$ such that
for initial data
$y^e_{+}$,  $y^e_{-} \in L^\infty(0,\, L^e)$ ($e\in E$)
such that
\[
\|y^e_{\pm}-J\|_{L^\infty(0, L^e)} \leq \varepsilon(T,\, M)
\]
and control functions
$u^e \in L^\infty(0,\, T)$  ($e\in E$)
such that
\[
\|u^e-J \|_{L^\infty(0, T)} \leq \varepsilon(T,\, M) 
\]
there exists a unique solution of
 ${\bf (S) }$ 
 that satisfies the integral equations
 (\ref{integralequation}) for all $e\in E$
 along the characteristic curves 
 with
 $S^e_{+}$,  $S^e_{-} \in L^\infty((0,\, L^e)\times (0,\, T))$ ($e\in E$)
 and the boundary condition
 (\ref{rbriemann2neu}) 
 at the boundary nodes 
 and the node condition
 (\ref{Omegav}) at the interior nodes 
 almost everywhere in $[0, \, T]$ 
such that for all $e\in E$ we have
\begin{equation}
\label{lunendlichbound}
\|S^e_{\pm} - J\|_{L^\infty( (0,\, T) \times (0,\, L^e) ) } \leq M.
\end{equation}
This solution depends in a stable way  on the initial and boundary data in the sense that
for
initial data
$
\|z^e_{\pm}-J\|_{L^\infty(0, L^e)} \leq \varepsilon(T,\, M)
$
and control functions
$v^e \in L^\infty(0,\, T)$  ($e\in E$)
such that
$
\|v^e-J \|_{L^\infty(0, T)} \leq \varepsilon(T,\, M)
$
for the corresponding solution $\tilde S^e_\pm$ we have the inequality
\[
\| \tilde S^e_\pm - S^e_{\pm}\|_{L^\infty( (0,\, T) \times (0,\, L^e) ) } \leq 
C(T)\; \max_{e \in E} \{ \|y^e_{\pm} - z^e_{\pm}\|_{L^\infty(0, L^e)}; \|u^e - v^e \|_{L^\infty(0, T)}\}
\]
where $C(T)> 0$ is a constant that does not depend on $z^e_{\pm}$ or $v^e$.
If 
\begin{equation}
    \label{kleinezeit}
T \leq \min_{e\in E} \frac{L^e}{c} 
\end{equation}
the solution satisfies the a  priori bound
\begin{equation}
    \label{aprioribound2020}
 \mbox{\textrm ess}
 \sup_{s \in [0,\, T]}
 \max_{e\in E}
\left\{
\|R^e_+(s,\, x) - J\|_{L^\infty(0, \, L^e)},
\;
\|R^e_-(s,\,x)- J \|_{L^\infty(0, \, L^e)}
\right\}
\end{equation}
\[
\leq
\varepsilon(T,\, M) \exp( 16  \, \max_{e\in E} \nu^e\, M \,T)
.\]

\end{theorem}
\textbf{Proof:} 
The proof is based upon
Banach's fixed point theorem
with the canonical fixed point iteration.
It has to be shown that this map is a contraction in the Banach space
$ X = \times_{e \in E} 
\left( L^\infty(  (0,\, T) \times (0,\, L^e)  ) \right)^2 
$
on the set $B(M)$ 
\[
= \times_{e\in E} 
\left\{ (S^e_{+},  \; S^e_{-})\in X:
\; (\ref{lunendlichbound}) \; \mbox{\rm 
and the initial conditions of
 ${\bf (S) }$ at $t=0$ hold}
\right\}.
\]
 In order to show this, we  use 
 an upper bound for the source term
 in  (\ref{integralequation})  that is 
 given by the
 continuously differentiable function $\sigma^e(R^e_+,\, R^e_-)$.
 In fact,  for $R\in B(M)$ 
 for all $e\in E$ due to 
 \eqref{sigmadefinition} we have
 $|\sigma^e (R^e_+,\, R^e_-)| \leq 4\,\nu^e \,
 M^2$. 
 
Moreover, it has to be shown that the iteration map maps
from $B(M)$ into $B(M)$. This is true if $M$
and $\varepsilon(T,\, M)$ are chosen sufficiently small.

In this analysis, it
has to be taken into account that
the characteristic curves can 
be reflected at the
boundaries of
the edges a finite number of times.
Due to  the linear node condition (\ref{Omegav}) 
in each such crossing 
the absolute value of
the outgoing  Riemann invariants
can be at most three times as large as 
the largest absolute values of the incoming  Riemann invariants.
For $t\in [0, T]$ almost everywhere 
and $v\in V$ 
with $|E_0(v)|\geq 2$ 
we have the inequality 
\begin{equation}
\label{knotenungleichung}
|S^e_{out}(t,\, x^e(v)) - J| \leq 
3 \, \max_{f \in E_0(v)} |S^f_{in}(t, \, x^e(v) ) - J|.
\end{equation}

In a first step, we assume that 
the time horizon is sufficiently small in the
sense that 
\begin{equation}
\label{Tvoraussetzung}
T < \frac{1}{16} \, \frac{1}{\max_{e\in E} \nu^e\, M}.
\end{equation}

holds. 
For given ${({\cal R}^e})_{e\in E}
= ((R^e_+,\, R^e_-))_{e\in E}\in B(M)$
we define
\begin{eqnarray*}
\Phi_+^e(R^e_+,\, R^e_-)(t,x)
& = &
\Xi^e_+(t^{e}_+(x,t),\,\xi^{e}_+(t^{e}_+(x,t),\, x,\, t))
\\
& - &
\int_{t^{e}_+(x,t)}^t  \sigma^e(R^e_+,\,R^e_-) \,\left(\tau,\, \xi^{e}_+(\tau,\, x, \, t) \right) \, d \tau
\end{eqnarray*}
with $\sigma^e$ as defined in
(\ref{sigmadefinition}).
Here we define
\[
\Xi^e_+(t^{e}_+(x,t),\,\xi^{e}_+(t^{e}_+(x,t),\, x,\, t))
\]
\[
=
\left\{
\begin{array}{r}
(1 -  \mu^v) u_+^e(t^{e}_+(x,t)) 
+ \mu^v R^e_-(t^{e}_+(x,t),\,\xi^{e}_+(t^{e}_+(x,t),\, x,\, t))
  \\
{\rm if}  \; t^{e}_+(x,t)>0,\, 0=x^e(v)\;{\rm  and }\; |E_0(v)|=1;
\\
y_+^e( \xi^{e}_+(0,\, x,\, t))
\;\; {\rm if} \;  t^{e}_+(x,t)=0;
\\
\Omega_v^e \,  R^v_{in} ( t^{e}_+(x,t) )
\;\;
{\rm if} \;  t^{e}_+(x,t)>0, \, 0=x^e(v)\; {\rm  and }\; |E_0(v)|\geq 2.
\end{array}
\right.
\]
Here $\Omega^v$ is the square
matrix that describes the linear 
interior node conditions
(\ref{Omegav}).

The components of
$R^v_{in} (t)$
that appear in the last line are in turn obtained by integrating
along the characteristic curves $\xi^{e}_\pm$. 
Due to (\ref{kleinezeit}), they
can be followed back to the initial state,
that is
for $f\in E_0(v)$
the components of $R^v_{in} (t)$ have the form
\begin{equation}
\label{strichrvinintegration}
 R^f_{\pm} (t,\, x^f(v))
%
 =
y_\pm^f( \xi^{f}_\pm(0,\, x,\, t))
 \mp
\int_{0}^t  \sigma^f(R^f_+,\,R^f_-) \,\left(\tau,\, \xi^{f}_\pm(\tau,\, x, \, t) \right) \, d \tau
\end{equation}
without further reflections.
%
Analogously we define
\begin{eqnarray*}
\Phi_-^e(R^e_+,\, R^e_-)(t,x)
& = &
\Xi^e_-(t^{e}_-(x,t),\,\xi^{e}_-(t^{e}_-(x,t),\, x,\, t))
\\
& + &
\int_{t^{e}_-(x,t)}^t  \sigma^e(R^e_+,\,R^e_-) \,\left(\tau,\, \xi^{e}_-(\tau,\, x, \, t) \right) \, d \tau.
\end{eqnarray*}
Here we define
\[
\Xi^e_-(t^{e}_+(x,t),\,\xi^{e}_-({e}_+(x,t),\, x,\, t))
\]
\[
=
\left\{
\begin{array}{r}
(1 - \mu^v) u_-(t^{e}_-(x,t)) 
+ \mu^v R^e_+(t^{e}_+(x,t),\,\xi^{e}_-({e}_+(x,t),\, x,\, t))
\\
{\rm if } \; t^{e}_-(x,t)>0, \, L^e=x^e(v)\;{\rm  and }\; |E_0(v)|=1;
\\
y_-( \xi^{e}_-(0,\, x,\, t))
\;\;{\rm if } \, t^{e}_-(x,t)=0; 
\\
\Omega_v^e \,  R^v_{in} (t^{e}_-(x,t)) 
\;\; 
{\rm if } \, t^{e}_-(x,t)>0, \, L^e=x^e(v)\;{\rm  and }\; |E_0(v)| \geq 2.
\end{array}
\right.
\]
Again the components of
$R^v_{in} (t)$ are obtained by integrating
along the corresponding characteristic curves $\xi^{f}_\pm$
for $f\in E_0(v)$
going back to the given initial values as in
(\ref{strichrvinintegration}).

In this way we get the fixed point iteration
where for all $e\in E$ we define
\begin{equation}
\left(
\begin{array}{c}
\rho^{e,(k+1)}_+(t,x)
\\
\rho^{e,(k+1)}_-(t,x)
\end{array}
\right)
=
\left(
\begin{array}{c}
\Phi_+^e( \rho^{e,(k)}_+,\, \rho^{e,(k)}_-)(t,x)
\\
\Phi_-^e( \rho^{e,(k)}_+,\, \rho^{e,(k)}_-)(t,x)
\end{array}
\right)
\end{equation}
and that we start with functions
$(\rho^{e, (1)}_+,\, \rho^{e, (1)}_-)_{e\in E}\in B(M)$.
%
Our aim is to apply Banach's fixed point theorem.
We check  in several steps that the assumptions hold.
First we show that the fixed point iteration is well--defined.

{\bf Step 1
(The fixed point iteration is well--defined)}
In order to show that the fixed point iteration is well--defined,
we show that the iterates  remain in $B(M)$.
Assume that  $\rho^{e,(k)}
=
\left(
\rho^{e,(k)}_+,\, \rho^{e,(k)}_-
\right)
\in  B(M)$.
%
%
%
%
For $e\in E$,  define
\begin{equation}
\label{supdefo}
 {\cal S}^{e,(k+1)}(t):=\mbox{\rm ess}\sup_{ x\in [0,L^e]}  \left\{|\rho^{e,(k+1)}_\pm(t,x)- J|
 \right\}
 \;
 \mbox{ \textrm{and} }
 {\cal S}^{(k+1)}(t):=\max_{e\in E} 
{\cal S}^{e,(k+1)}(t).
\end{equation}
As long as there is
at most one crossing
of a characteristic curve through an edge 
the definition of {\bf (S)} implies 
\[
  {\cal S}^{e,(k+1)}(t) \leq 
3 \, 
\varepsilon(T,\, M)
+
16  \, \max_{e\in E} \nu^e\, M^2\, T.
%
\]
Define $\bar \nu = \max_{e\in E} \nu^e$. 
Then we have
 $
{\cal S}^{e,(k+1)}(t)
\leq
3 \, 
\varepsilon(T,\, M)
+
16  \, \bar \nu\, M^2\, T
%
$. 
Thus we have
\[
\mbox{\rm ess}\sup_{t\in [0,T], \, x\in [0,L^e]} |\rho^{e,(k+1)}_\pm(t,x)- J|
\leq
3 \, 
\varepsilon(T,\, M)
+
16  \, \bar \nu\, M^2\, T.
\]
Now $M$ and $\varepsilon(T,\, M)$
have to be chosen in  such a way that 
\begin{equation}
3 \, 
\varepsilon(T,\, M)
+
16  \, \bar \nu\, M^2\, T
    \leq M.
\end{equation}
Due to 
\eqref{Tvoraussetzung} 
this is possible for 
$\varepsilon(T,\, M)
=
\frac{M}{3} 
\left( 1 - 16  \, \bar \nu\, M\, T \right)>0
$. 
Then we have
\[|\rho^{e,(k+1)}_\pm(t,x)- J|\leq M.\]
By induction this implies that for all $k\in \{1,2,3,...\}$
we have
$( \rho^{e,(k)}_+,\, \rho^{e,(k)}_-)\in B(M)$. 
Hence 
all the iterates of 
the fixed point iteration 
remain in the set $B(M)$.


%
%
%

{\bf Step 2:  Contractivity}
The next step is to show that $\Phi$ is a contraction.
Let
$(R^e_+,\,R^e_-)$, $(S^e_+,\,S^e_-)\in 
B(M)$.
For $(t,x)\in [0,T]\times [0, L^e]$,
the definition of $\Phi_{\pm}$ implies the inequality
\begin{eqnarray*}
\left|
\Phi^e_{\pm}(R^e_+,R^e_-) - \Phi^e_{\pm}(S^e_+,S^e_-)
\right|(t,x)
\leq A^e + I^e,
\end{eqnarray*}
with
\begin{eqnarray*}
A^e
& = &
\left|
R^e_\pm(t^e_\pm(x,t),\,\xi^{{e}}_\pm(t^{{e}}_\pm(x,t),\, x,\, t))
-
S^e_\pm(t^e_\pm(x,t),\,\xi^{{e}}_\pm(t^{e}_\pm(x,t),\, x,\, t))
\right|,
\\
I^e
&
=
&
\int_0^t
\left|
 \sigma^e(R^e_+,\,R^e_-)\left(\tau,\, \xi^{e}_\pm(\tau,\, x, \, t) \right)
 -
 \sigma^e(S^e_+,\,S^e_-)
  \left( \tau,\, \xi^{{e}}_\pm(\tau,\, x, \, t) \right)
 \right|\, d \tau.
\end{eqnarray*}
We have the inequality
\begin{eqnarray*}
I^e
 &
\leq
&
\int_0^t\, 4 \,\nu^e \, M  \,
\left|
 R^e_+\left(\tau,\, \xi^{{e}}_\pm(\tau,\, x, \, t) \right)
 -
 S^e_+\left(\tau,\, \xi^{{e}}_\pm(\tau,\, x, \, t) \right)
 \right|
 \, d \tau
 \\
 &
 +
&
 \int_0^t\,  4 \,\nu^e \, M  \,
\left|
 R^e_-\left(\tau,\, \xi^{{e}}_\pm(\tau,\, x, \, t) \right)
 -
 S^e_-\left(\tau,\, \xi^{{e}}_\pm(\tau,\, x, \, t) \right)
 \right|
  \, d \tau
 \\
 &
\leq
&
4 \,T\, \nu^e \, M  \,
 \left\|
{{\cal R}^e}
 -
 {{\cal S}^e}
 \right\|_{L^\infty([0,T]\times [0,L^e])^2}
.
\end{eqnarray*}
Hence we have  the inequality
\begin{eqnarray*}
I^e
&
\leq
&
4 \,T\, \nu^e \, M  \,
\|{{\cal R}^e}- {{\cal S}^e} \|_{L^\infty([0,T]\times [0,L^e])^2}.
\end{eqnarray*}

%
Now we look at the term $A^e$.
We consider three cases
i)-iii). 
{\em Case i): }
If $t^e_{\pm}(x,\, t)=0$, 
due to the initial conditions,
we have $A^e=0$.
{\em Case ii): }
If $t^e_{\pm}(x,\, t)>0$,
$P^{e\pm}_0(t,\, x)=
(t^e_{\pm}(x,\, t), \, x^e(v))$
and for the corresponding $v\in V$
we have  $|E_0(v)|=1$, then 
due to the fact that 
(\ref{kleinezeit}) implies that 
the characteristics
entering $v$ can be traced back
directly to the initial time,
and $|\mu^e|<1$
we have 
\[
|A^e|
\leq |\mu^e| \max_{f\in E} \; I^f
\leq 12 \, T \, M \,  \max_{f\in E} \; \nu^f
\|{{\cal R}^f}- {{\cal S}^f} \|_{L^\infty([0,T]\times [0,L^f])^2}.
\]
{\em Case iii): } 
Similarly, if $t^e_{\pm}(x,\, t)>0$,
$P^{e\pm}_0(t,\, x)=
(t^e_{\pm}(x,\, t), \, x^e(v))$
and for the corresponding $v\in V$
we have  $|E_0(v)|>1$,
due to the fact that 
(\ref{kleinezeit}) implies that 
there is at most one crossing through
a node that can result at most by 
an increase by the factor $3$, 
we obtain
$|A^e|
\leq
12 \, T \,  M \,  \max_{f\in E} \; \nu^f
\|{{\cal R}^f}- {{\cal S}^f} \|_{L^\infty([0,T]\times [0,L^f])^2}$.

Hence for the term $A^e$
we have the Lipschitz constant $12\, \max_{e\in E}
T\, \nu^e \, M $.


With our results for $I^e$ and $A^e$ we obtain the Lipschitz inequality for $\Phi^e_\pm$
\begin{eqnarray*}
& &
\left\|
\Phi^e_{\pm}(R^e_+,\, R^e_-) - \Phi^e_{\pm}(S^e_+,S^e_-)
\right\|_{L^\infty([0,\,T]\times [0,L^e])}
\\
&
\leq
&
 I^e  + A^e
 \leq 
\max_{f\in E}
\left[
16 \,T\, \nu^f \, M  
\right]
\,
 \|{{\cal R}^f}- {{\cal S}^f} \|_{L^\infty([0,T]\times [0,L^f])^2}
.
\end{eqnarray*}

With the notation $Y^e = L^\infty([0,T]\times [0,L^e])$ this implies
\[
\max_{e\in E}
\left\{
\left\|
\Phi^e_{+}(R^e_+,R^e_-) - \Phi^e_{+}(S^e_+,S^e_-)
\right\|_{Y^e},
\;
\left\|
\Phi^e_{-}(R^e_+,R^e_-) - \Phi^e_{-}(S^e_+,S^e_-)
\right\|_{Y^e}
\right\}
\]
\begin{eqnarray*}
&
\leq &
L_{\text kontr}\;
\max_{e\in E}
\left\{
 \left\|
R^e_+ - S^e_+
\right\|_{Y^e},
\;
 \left\|
R^e_-  - S^e_-
\right\|_{Y^e}
\right\}
\end{eqnarray*}
with the contraction constant
$
L_{\text kontr}
=
\max_{e\in E}
\left[
16 \,T\, \nu^e \, M
\right].
$
Due to \eqref{Tvoraussetzung}
we have
$L_{\text kontr}<1$.
%
%
Hence the map
$\Phi=(\Phi_+,\Phi_-)$
is a contraction. 
Thus Banach's fixed point theorem implies
the existence of a unique fixed point of the map,
which solves our semilinear initial boundary value problem  ${\bf (S)}$
if $T$ satisfies 
\eqref{Tvoraussetzung} and 
(\ref{kleinezeit}).
For $t\in [0,\, T]$,  define the number
\[
 U(t) =
 \mbox{\textrm ess}
 \sup_{s \in [0,\, t]}
 \max_{e\in E}
\left\{
\|R^e_+(s,\, x) - J\|_{L^\infty(0, \, L^e)},
\;
\|R^e_-(s,\,x)- J 
\|_{L^\infty(0, \, L^e)}
\right\}
 .
 \]
Since the solution $R^e_\pm$ is a fixed point
of $\Phi^e_\pm$, the definition of $\Phi^e_\pm$
 implies the integral inequality
\[U(t) \leq  U(0) +
\int_0^t
16  \, \max_{e\in E} \nu^e\, M \, 
\, U(\tau)
\,d\tau
\]
for all $t\in [0,\, T]$.
Now we can apply Gronwall's Lemma (see for example \cite{gubook})
and obtain for all $s\in [0,\, T]$
the upper bound
\[
U(s) \leq
U(0)\;
\exp( 16  \, \max_{e\in E} \nu^e\, M \,s)
\leq
\varepsilon(T,\, M) \exp( 16  \, \max_{e\in E} \nu^e\, M \,T)
.
\]
Thus we have shown
the a--priori bound 
(\ref{aprioribound2020})
for sufficiently small time-horizons
$T$ that satisfy 
\eqref{kleinezeit} and \eqref{Tvoraussetzung}.

For arbitrarily large $T>0$,
we obtain a solution in the following way:
Define $T_0 = \min_{e\in E} L^e/c.$
Then we obtain a solution on
the interval $[0,\, T_0]$
as shown  above,
and  the a priori bound
(\ref{aprioribound2020})
yields the bound
\begin{equation}
    \label{bound2020}
\|R^e_{\pm}(T_0, \cdot) -J\|_{L^\infty(0, L^e)} \leq 
\varepsilon(T,\, M) \exp( 16  \, \max_{e\in E} \nu^e\, M \,T_0)
\end{equation}
for the state  at  time $T_0$.
Define 
$
\varepsilon_1(T,\, M) 
=
\varepsilon(T,\, M) \exp( -16  \, \max_{e\in E} \nu^e\, M \,T_0).
$
Then if we start with data that satisfy 
\[
\|y^e_{\pm}-J\|_{L^\infty(0, L^e)} \leq \varepsilon_1(T,\, M)
\]
and control functions
$u^e \in L^\infty(0,\, 2 \,T_0)$  ($e\in E$)
such that
\[
\|u^e-J \|_{L^\infty(0, 2\,T_0)} \leq \varepsilon_1(T,\, M) 
\]
we obtain first a solution
on $[0, \, T_0]$.
Due to 
(\ref{bound2020})
and the definition of 
$\varepsilon_1(T,\, M) $
we can use the same argument
again to obtain the solution on
the time interval
$[T_0,\, 2\, T_0]$.
More generally,  
for $n \in \{1,2,3,...\}$ if the data satisfy
\[
\|y^e_{\pm}-J\|_{L^\infty(0, L^e)} \leq \varepsilon_n(M)\;
\mbox{  \textrm{and}  }  
\|u^e-J \|_{L^\infty(0,\, n\, T_0)} \leq \varepsilon_n(T,\, M) 
\]
where
$
\varepsilon_n(T,\, M) 
=
\varepsilon(T,\, M) \exp( -16 n \, 
\max_{e\in E} \nu^e
\, M \,T_0)
$ then
we obtain a solution
on the interval
$[0,\, n \, T_0]$.
%
%
Thus we have proved
Theorem \ref{existence}. 
$\Box$


\begin{remark}
\label{existenceremark}
An analogous existence result holds for 
more regular solutions
in $C([0,\, T], \, H^1(0, L^e)) $, ($e\in E$)
(in the space $C([0,\, T], \, W^{1,\infty}(0, L^e)) $, ($e\in E$)
respectively)
 for initial and boundary data
  that are compatible with each other and
with the node conditions
 such that $y^e_{\pm}-J$ and  $u^e-J$
 have  sufficiently small
 norms in $H^1$ ($W^{1,\infty}$ respectively).
 
The existence  result for solutions with $H^2$-regularity
(more precisely in the space 
$\times_{e\in E} C([0,\, T], \, H^{2}(0, L^e)) $)
for quasilinear systems
given in \cite{bastincoron}
requires $C^2$-regularity of the source term.
In our case, the source term only has $C^1$-regularity.
Therefore, the result from \cite{bastincoron} cannot be applied.
\end{remark}



For the proof of the exponential decay
of the $L^2$-norm of $\delta$,
we need an observability inequality
for the $L^2$-norm which 
is presented in Section \ref{l2observability}.
An observability inequality
for the $H^1$-norm is shown in
Section \ref{observabilityh1}.
It allows
to analyze
the  exponential decay
of the $H^1$-norm of $\delta$,


\section{An $L^2$--observability inequality for a  network}
\label{l2observability}

In this section we derive an observability inequality for a  network.
The aim is to get an upper bound for the $L^2$-norm of the system
state on any edge $e$ at the time $t$ in terms of the
$L^2$-norm of
the trace of the state at one of the   nodes adjacent to $e$
on the time interval $[t- T,\, t + T]$ with $T>0$ sufficiently large.
In \cite{zuazua1}, an observability inequality for a star-shaped network of
strings (without source term) is derived.



\begin{theorem}
\label{observabilitythm}
Let $v \in V$, $e \in E_0(v)$.
Assume that $T>\frac{L^e}{c}$ and $t>T$.
Assume that system {\bf (Diff)} has a solution on $[0,\, t + T]$
such that for  $x$ almost everywhere in $[0, \, L^e]$ we have
$\delta^e_+(\cdot, \, x)$, $\delta^e_-(\cdot, \, x)\in L^2(0,\, T + t)$
and that there exists a constant
$\tilde M $ such that
for almost all
 $x$ in $[0, \, L^e]$ 
for the solution of {\bf (S)} we have the inequalities
\begin{equation}
\label{druckvor}
|S^e_+(s,\, x) - S^e_-(s,\, x)| \leq \tilde M, \;\;
|S^e_+(s,\, x) - S^e_-(s,\, x) + \delta^e_+(s,\, x) - \delta^e_-(s,\, x)| \leq \tilde M
\end{equation}
for $s$  almost everywhere in   $[0,\, t + T]$.
Then,
there exists a constant $C_0^e(\tilde M)$ such that
for all $t >T$
the following inequality holds:
\begin{multline}
\label{observabilityinequality}
 \int_0^{L^e}
|\delta^e_+(t,\, x)|^2 + |\delta^e_-(t,\, x)|^2 \, dx\\
\leq
C_0^e(\tilde M) \, \int_{t- T}^{ t + T}
|\delta^e_+(s,\, x^e(v) )  |^2 + |\delta^e_-(s,\, x^e(v) ) |^2 \, ds
\end{multline}

\end{theorem}

\begin{remark}
In Theorem \ref{observabilitythm}, we observe from one end $x^e(v)$ of $e$ only. If we observe from both sides, the constant improves by a factor of $2$.
\end{remark}

\begin{remark}
Since (\ref{lunendlichbound}) with
$M = \frac{1}{2} \, \tilde M$ implies
(\ref{druckvor}) ,
Theorem \ref{existence}
applied to $R$ and $S$ yields sufficient conditions
for (\ref{druckvor}) if $\tilde M$ is sufficiently small. 

An a priori upper bound
(\ref{druckvor})  also holds 
for  classical solutions
(even in the sense of a maximum), see \cite{libook}.
Also for solutions in 
$C([0, \, T], \; H^1(0,\, L^e))$, 
(\ref{druckvor})  holds.
\end{remark}

\textbf{Proof:}
Let us give the proof of \eqref{observabilityinequality}: 

For $v \in V$, $e\in E_0(v)$, $t \geq T$ and $x\in [0, L^e]$
consider

\begin{equation}
\label{hdefinition}
{\cal H}^e(x)
=
\frac{1}{2}
\int_{t- \frac{L^e - |x^e(v) -  x| }{c}}^
{t + \frac{L^e - |x^e(v) -  x|}{c}}
\left|\delta^e_+(s,x) \right|^2 + \left|\delta^e_-(s,x)\right|^2
\, ds \geq 0.
\end{equation}

Then we have  
${\cal H}^e(L^e - x^e(v)  ) = 0$.

In order to make
the proof more readable, we only give it in case $x^e(v)= L^e$.
The case $x^e(v)= 0$ is analogous.
For the derivative of ${\cal H}^e$ with respect to $x$ we have almost everywhere 
\begin{multline*}
 \frac{d}{dx}
{\cal H}^e(x)
 =
\int_{t- \frac{x}{c}}^{t + \frac{x}{c}}
\delta^e_+(s,x)\, \partial_x \delta^e_+(s,x) + \delta^e_-(s,x) \,\partial_x  \delta^e_-(s,x)
\, ds
\\
 +   \frac{1}{2\, c}\left[ \left|\delta^e_+(t + \frac{x}{c} \, ,x)\right|^2 
+ \left|\delta^e_-(t + \frac{x}{c} \, ,x) \right|^2
+
\left|\delta^e_+(t - \frac{x}{c} \, ,x)\right|^2 + \left|\delta^e_-(t - \frac{x}{c} \, ,x) \right|^2
\right].
\end{multline*}
Due to the partial differential equation in system ${\bf (Diff)}$   this yields
$\frac{d}{dx}{\cal H}^e(x)$
\begin{eqnarray}
& = & \nonumber
\int_{t- \frac{x}{c}}^{t + \frac{x}{c}}
- \delta^e_+(s,x)\, \frac{1}{c} \partial_t \delta^e_+(s,x) 
+  \frac{1}{c} \, \delta^e_-(s,x)  \,\partial_t  \delta^e_-(s,x)
\\
&&  \nonumber
\quad -\frac{1}{c} \, (\delta^e_+  + \delta^e_- )       \, 
\left[
\sigma^e (\delta_+^e+S_+^e,\,\delta_-^e + S_-^e)
-
\sigma^e (S_+^e,\,  S_-^e)
\right]
\, ds
\\
& & \quad \nonumber +\frac{1}{2\, c}\left[ \left|\delta^e_+(t + \frac{x}{c} \, ,x)\right|^2 + \left|\delta^e_-(t + \frac{x}{c} \, ,x)  \right|^2
+
\left|\delta^e_+(t - \frac{x}{c} \, ,x)\right|^2 +
 \left|\delta^e_-(t - \frac{x}{c} \, ,x)  \right|^2
\right]
\\
& = &\nonumber
\int_{t- \frac{x}{c}}^{t + \frac{x}{c}}
-  \frac{1}{2c} \partial_t \left|\delta^e_+(s,x)\right|^2  
+  \frac{1}{2c} \partial_t \left|\delta^e_-(s,x) \right|^2
\\
&& \quad \nonumber -\frac{1}{c} \, (\delta^e_+  + \delta^e_- )       \, 
\left[
\sigma^e (\delta_+^e+S_+^e,\,\delta_-^e + S_-^e)
-
\sigma^e (S_+^e,\,  S_-^e)
\right]
\, ds
\\ \nonumber
&  &  \quad +\frac{1}{2\, c}\left[ \left|\delta^e_+(t + \frac{x}{c} \, ,x)\right|^2 + \left|\delta^e_-(t + \frac{x}{c} \, ,x)  \right|^2
+
\left|\delta^e_+(t - \frac{x}{c} \, ,x)\right|^2 +
\left|\delta^e_-(t - \frac{x}{c} \, ,x) \right|^2
\right]
\\ \nonumber
& = &
  \frac{1}{2c} \left[  - \left|  \delta^e_+(s,x) \right|^2
  +   \left|\delta^e_-(s,x)\right|^2 \right]\Bigg|_{s=t- \frac{x}{c}}^{t + \frac{x}{c}}
\\
&& \quad \nonumber -  \int_{t- \frac{x}{c}}^{t + \frac{x}{c}}\frac{1}{c} \, (\delta^e_+ 
+ \delta^e_- )       \, \left[
\sigma^e (\delta_+^e+S_+^e,\,\delta_-^e + S_-^e)
-
\sigma^e (S_+^e,\,  S_-^e)
\right]
\, ds
\\  \nonumber
&  & \quad + \frac{1}{2\, c}\left[ \left|\delta^e_+(t + \frac{x}{c} \, ,x)\right|^2 
+ \left|\delta^e_-(t + \frac{x}{c} \, ,x) \right|^2
+
\left|\delta^e_+(t - \frac{x}{c} \, ,x)\right|^2 
+ \left|\delta^e_-(t - \frac{x}{c} \, ,x) \right|^2
\right]
\\ \nonumber
& \color{black}{\geq} &
{-} \int_{t- \frac{x}{c}}^{t + \frac{x}{c}}\frac{1}{c} \, \left[   
 2 \, \nu^e \, \tilde M 
 \right]     
{\color{black}{|}}\delta^e_+  + \delta^e_- {\color{black}{|}} \, 
{\color{black}{|}}\delta^e_+ -  \delta^e_-{\color{black}{|}}
\, ds
\\ 
\nonumber
&  & 
    + \left[\frac{1}{2\, c}  
    +  \frac{1}{2c}\right] \left|\delta^e_-( t + \frac{x}{c},x)  \right|^2
   +  \left[\frac{1}{2\, c}  +   \frac{1}{2c}\right]
    \left|\delta^e_+(t - \frac{x}{c}  ,x) \right|^2 
    .
\end{eqnarray}
Here the last inequality follows with the mean value theorem
from the definition 
(\ref{sigmadefinition}) 
of the function $\sigma^e$,
since the function $z\mapsto \nu^e \, z\, |z|$ is
differentiable with the derivative
$z \mapsto 2\, \nu^e \, |z|$.
Then, assumption (\ref{druckvor}) 
is applied to obtain an upper bound.
This implies the inequality
\begin{equation}
\label{eq:deltaboundary}
\frac{d}{dx} {\cal H}^e(x)
\geq
-  4\, \frac{\nu^e}{c}  \, \tilde M \,
{\cal H}^e(x)
\end{equation}
for $x\in [0,\,L^e]$ almost everywhere.
%
%
%
With (\ref{druckvor}),
due to Gronwall's Lemma,
this implies
the inequality
\begin{equation}
\label{2201h}
{\cal H}^e(x) \leq \exp \left(4\, \frac{\nu^e}{c} \tilde M \, (L^e  - x) \right) \,  {\cal H}^e(L^e)
\end{equation}
for all $x\in [0,\, L^e]$.
Since for real numbers $r_1$, $r_2$ we have
$(r_1 + r_2)^2 \leq 2 \, r_1^2 + 2 \, r_2^2$,
for  
$x\in [0,\, L^e]$ almost everywhere 
 we obtain 
\begin{eqnarray*}
& & |\delta^e_+(t,\, x)|^2
\\
&
=
&
 \left|\delta^e_+( t +  \frac{L^e - x}{c}, L^e) 
 - \frac{1}{c}
   \int_{x}^{L^e}  \,\left(
\sigma^e (\delta_+^e+S_+^e,\,\delta_-^e + S_-^e)
-
\sigma^e (S_+^e,\,  S_-^e)
\right)
 ( t +  \frac{s - x}{c}, \, s)        \, ds \right|^2
\\
&
\leq
&
2 \, \left|\delta^e_+( t +  \frac{L^e - x}{c}, L^e) \right|^2 +
2  \, \left|  \frac{1}{c} \int_{x}^{L^e} 
\,\left(
\sigma^e (\delta_+^e+S_+^e,\,\delta_-^e + S_-^e)
-
\sigma^e (S_+^e,\,  S_-^e)
\right)
 ( t +  \frac{s - x}{c}, \, s)        \, ds \right|^2
 \\
 &
 \leq
 &
 2 \, \left|\delta^e_+( t +  \frac{L^e - x}{c}, L^e)\right|^2 +
 2  \, \left| \frac{1}{c}  \int_{x}^{L^e}  \,
2 \,\nu^e \, \tilde M \left[ |\delta^e_+| +  |\delta^e_-| \right]
 ( t +  \frac{s - x}{c}, \, s)        \, ds \right|^2
 \\
 &
 \leq
 &
 2 \, \left|\delta^e_+( t +  \frac{L^e - x}{c}, L^e)\right|^2 +
 16 \, \frac{1}{c^2} \, L^e  \, ( \nu^e )^2 \, (\tilde M)^2 \, \,   \int_{x}^{L^e}  \,
  \left[ |\delta^e_+|^2 +  |\delta^e_- |^2 \right]
 ( t +  \frac{s - x}{c}, \, s)        \, ds.
\end{eqnarray*}
Thus we have

\begin{eqnarray*}
& & \int_0^{L^e} |\delta^e_+(t,\, x)|^2 \, dx
\\
&
\leq
&
2 \, \int_0^{L^e}  \left|\delta^e_+( t +  \frac{L^e - x}{c}, L^e)\right|^2   \, dx
\\
&
+
&
 \frac{16}{c^2} \, L^e  \, ( \nu^e )^2 \, (\tilde M)^2
\int_0^{L^e}
 \int_{x}^{L^e}  \,
\left[ |\delta^e_+|^2 +  |\delta^e_- |^2 \right]
 ( t +  \frac{s - x}{c}, \, s)        \, ds
\, dx
\\
&
\leq
&
2 \, \int_0^{L^e}  \left|\delta^e_+( t +  \frac{L^e - x}{c}, L^e)\right|^2   \, dx
\\
&
+
&
\frac{16}{c} \, L^e  \, ( \nu^e )^2 \, (\tilde M)^2
\int_0^{L^e}
\int_{t- \frac{x}{c}}^{t + \frac{x}{c}}
\left|\delta^e_+(s,x) \right|^2 + \left|\delta^e_-(s,x)\right|^2
\, ds
\,dx.
\end{eqnarray*}
Hence,   using the definition (\ref{hdefinition}) of  ${\cal H}^e(x)$
and inequality (\ref{2201h}),  we obtain 
\begin{eqnarray*}
&& \int_0^{L^e} |\delta^e_+(t,\, x)|^2 \, dx
\\
&
\leq
&
2 \, \int_0^{L^e}  \left|\delta^e_+( t +  \frac{L^e - x}{c}, L^e) \right|^2   \, dx
+
\frac{32\, L^e  \, ( \nu^e \, \tilde M)^2 }{c}   
\int_0^{L^e}
 {\cal H}^e(x) \, dx \\
&
\leq
&
2 \, \int_0^{L^e}  \left|\delta^e_+( t +  \frac{L^e - x}{c}, L^e) \right|^2   \, dx
+
\frac{32}{c} \, L^e  \, ( \nu^e )^2 \, (\tilde M)^2
\int_0^{L^e}  \exp \left(4 \, \frac{\nu^e   \, \tilde M \, (L^e- x)}{c}  \right) \,dx \,  {\cal H}^e(L^e)
\\
&
\leq
&
2 \, c  \, \int_t^{t + \frac{L^e}{c}}  \left|\delta^e_+(s, \, L^e) \right|^2   \, ds
\\
&
&
+ 4 \,  L^e  \, \nu^e \, \tilde M \,
\exp  \left(4\, \frac{\nu^e}{c} \tilde M \, L^e \right) \, 
\int_{t - \frac{L^e}{c}} ^{t + \frac{L^e}{c}}  \left|\delta^e_+(s, \, L^e)\right|^2  +  
\left|\delta^e_-(s, \, L^e)\right|^2 \, ds
\\
& \leq &
\left[  2 \, c  +   4\, L^e  \, \nu^e \, \tilde M \,
\exp  \left(4 \, \frac{\nu^e   \, \tilde M \, L^e}{c}  \right)   \right] \, 
\int_{t - \frac{L^e}{c}} ^{t + \frac{L^e}{c}}  \left|\delta^e_+(s, \, L^e) \right|^2  +  \left|\delta^e_-(s, \, L^e) \right|^2 \, ds.
\end{eqnarray*}
Similarly as above,
this yields
\begin{multline*}
\int_0^{L^e} |\delta^e_-(t,\, x)|^2 \, dx
\\
\leq
\left[  2 \, c + 4\,  L^e  \, \nu^e \, \tilde M \,
\exp  \left(4\, \frac{\nu^e}{c} \tilde M \, L^e \right)   \right] \, 
\int_{t - \frac{L^e}{c}} ^{t + \frac{L^e}{c}}  \left|\delta^e_+(s, \, L^e) \right|^2  +  \left|\delta^e_-(s, \, L^e) \right|^2 \, ds.
\end{multline*}


Adding up the inequalities for $\delta^e_+$ and $\delta^e_-$ yields
the observability inequality  (\ref{observabilityinequality})
with
\begin{equation}
\label{c0definition}
C_0^e(\tilde M) =    2 \, \left[  2 \, c  + 4\,  L^e  \, \nu^e \, \tilde M \,
\exp  \left(4 \, \frac{ L^e  \, \nu^e  \, \tilde M }{c} \right)   \right] .
\end{equation}


\section{An $H^1$--norm--observability  inequality}
\label{observabilityh1}
Now we prove an observability inequality
where the norm of the time derivative is included. This  yields 
an observability inequality for the $H^1$-norm.

\begin{theorem}
\label{observabilitythmh1}
Assume that $T>\max_{e\in E} \frac{L^e}{c}$ and $t>T$.
Assume that systems  {\bf (Diff)} and {\bf (S)}   have a solution on $[0,\, t + T]$
such that for all $e\in E$ and $x$ almost everywhere in $[0, \, L^e]$ we have
$\delta^e_+(\cdot, \, x)$, $\delta^e_-(\cdot, \, x)$,
$S^e_+(\cdot, \, x)$, $S^e_-(\cdot, \, x)$,
$\partial_t \delta^e_+(\cdot, \, x)$, $\partial_t \delta^e_-(\cdot, \, x)$,
$\partial_t S^e_+(\cdot, \, x)$, $\partial_t S^e_-(\cdot, \, x)\in L^2(0,\, T + t)$
and that there exists a constant
$\tilde M $ such that \eqref{druckvor} holds
for all
$e\in E$ and 
for $s$  almost everywhere in   $[0,\, t + T]$
for the solutions of  {\bf (Diff)} and {\bf (S)}.
Assume that there exists
a real number $\tilde B>0$ such that 
\begin{equation}
    \label{janbedingung}
\|\partial_t (S^e_+ -
 S^e_-)(\cdot, \, x) \|_{L^\infty(
t- \frac{L^e - |x^e(v) -x| 
 }{c},\,
t+ \frac{L^e + |x^e(v) -x|
 }{c}
)}
\leq \tilde B
\end{equation}
for $x$ almost everywhere in $(0, \, L^e)$.
Then there exists a constant $C_1(\tilde M, \tilde B)$ such that
for all $v\in  V$, $e\in E_0(v)$ and 
for all $t >T$, we have the inequalities
\begin{multline}
\label{observabilityinequalityh1}
\int_0^{L^e}
|\partial_t \delta^e_+(t,\, x)|^2 + |\partial_t \delta^e_-(t,\, x)|^2 \, dx
+ 
 \int_0^{L^e}
| \delta^e_+(t,\, x)|^2 + | \delta^e_-(t,\, x)|^2 \, dx  \\
\leq
C_1(\tilde M,\, \tilde B)  \int_{t- T}^{ t + T}
|\partial_t \delta^e_+(s,\, x^e(v))|^2 + | \partial_t \delta^e_-(s,\, x^e(v)) |^2 \, ds
\\
+
C_1(\tilde M,\, \tilde B) 
 \int_{t- T}^{ t + T}
|\delta^e_+(s,\, x^e(v))|^2 +
| \delta^e_-(s,\, x^e(v)) |^2 \, ds.
\end{multline}
\end{theorem}

\begin{remark}\label{janremark1_observ}
Note that 
\eqref{janbedingung} assumes more regularity of the observed solution than  we can 
expect, i.e. than is guaranteed by our well-posedness result Theorem \ref{existence}.
It is possible to obtain the same result (with a larger constant $C_1$) if \eqref{janbedingung} is replaced by the 
weaker 
assumption that there exist $\epsilon>0$ and $\overline{B}<\infty$ such that
\begin{equation}\label{janbesser} \| \partial_t (S^e_+ - S^e_-) \|_{L^{2+\epsilon}([0,T] \times [0,L^e])}  \leq \overline{B}. \end{equation}
Note that \eqref{janbesser} is weaker than \eqref{janbedingung} but also not guaranteed by
Theorem \ref{existence}.

The key ingredient in the proof that uses 
\eqref{janbesser} 
is the observation that whenever $\partial_t (S^e_+ - S^e_-)$ appears it is multiplied by $\delta^e_+  - \delta^e_-$ and the fact that (in two space dimensions), every $L^p$-norm with $p< \infty$ is controlled by the $H^1$-norm.

In what follows, we present the proof under the stronger assumption \eqref{janbedingung}
 as it is clearer and seems more appropriate to convey the main ideas.\end{remark}

{\bf Proof.}
In order to prove \eqref{observabilityinequalityh1}, we note  that the evolution of the space derivatives 
$(\delta^e_\pm)_x$
is governed by the following
partial differential equations:
\begin{multline}\label{eq:dxdeltaevol}
\partial_{x t}\left( \delta^e_+ \right)
=
\frac{1}{c}
\left[
-
\partial_{t t}\left( \delta^e_+ \right)
-
2\, \nu^e
\,|\delta^e_+ - \delta^e_- + S^e_+ - S^e_-|
\,
\partial_t (\delta^e_+ - \delta^e_- + S^e_+ - S^e_-) \right.\\
\left.+
2\, \nu^e\,
|S^e_+ - S^e_- |  \, 
 \partial_t (S^e_+ - S^e_-)
\right]  
\end{multline}
\begin{multline}
\partial_{x t}\left( \delta^e_- \right)_
=
\frac{1}{c}
\left[
\partial_{t t}\left( \delta^e_+ \right)
-
2\, \nu^e
\,|\delta^e_+ - \delta^e_- + S^e_+ - S^e_-|
\,
\partial_{t}(\delta^e_+ - \delta^e_- + S^e_+ - S^e_-)
\right.\\
\left.+
2\, \nu^e\,
|S^e_+ - S^e_- |  \, 
\partial_{t}(S^e_+ - S^e_-) 
\right].
\end{multline}
For $v\in V,\, e\in E_0(v)$, $t \geq T$ and $x\in [0, L^e]$
consider
\begin{equation}
\label{kdefinition}
{\cal K}^e(x)
:=
\frac{1}{2}
\int_{t- \frac{L^e - |x^e(v) -x| 
 }{c}}^{t+ \frac{L^e - | x^e(v) -x|
 }{c}}
\left|\partial_t \delta^e_+(s,x)\right|^2 + \left| \partial_t \delta^e_-(s,x)\right|^2
\, ds \geq 0.
\end{equation}
Then we have  ${\cal K}^e( L^e - x^e(v)  ) = 0$.

Now for the sake of
readability, we consider the case
$x^e(v)=L^e$.
For the derivative of ${\cal K}^e$ with respect to $x$ we have almost everywhere
\begin{eqnarray*}
& & \frac{d}{dx}
{\cal K}^e(x)
 =
\int_{t- \frac{x}{c}}^{t + \frac{x}{c}}
\partial_t \, \delta^e_+(s,x)\, \partial_{x t} \delta^e_+(s,x) + \partial_t \, \delta^e_-(s,x)\,\partial_{x t} \, \delta^e_-(s,x)
\, ds
\\
& + &  \frac{1}{2\, c}\left[\left|\partial_t \, \delta^e_+(t + \frac{x}{c} \, ,x)\right|^2 + \left| \partial_t \, \delta^e_-(t + \frac{x}{c} \, ,x) \right|^2
+
\left|\partial_t \, \delta^e_+(t - \frac{x}{c} \, ,x) \right|^2 + \left|\partial_t \, \delta^e_-(t - \frac{x}{c} \, ,x) \right|^2
\right].
\end{eqnarray*}
Due to \eqref{eq:dxdeltaevol}  this yields
$\frac{d}{dx}{\cal K}^e(x)$
\begin{eqnarray}  \nonumber
& = & \nonumber
\int_{t- \frac{x}{c}}^{t + \frac{x}{c}}
- \frac{1}{c} \, ( \partial_t \, \delta^e_+(s,x) )\, \partial_{ t t} \delta^e_+(s,x)
  +  \frac{1}{c} (\partial_t \, \delta^e_-(s,x)) \, \partial_{ t t}   \delta^e_-(s,x)
\\
& & \nonumber
-\frac{1}{c} \, ( \partial_t \, \delta^e_+ + \partial_t \, \delta^e_- )  \,
\partial_t\left(
\sigma^e (\delta_+^e+S_+^e,\,\delta_-^e + S_-^e)
-
\sigma^e (S_+^e,\,  S_-^e)
\right)
\,
\, ds
\\  \label{untilheresame}
&  &  
\nonumber
+\frac{1}{2\, c}\left[ \left| \partial_t \, \delta^e_+(t + \frac{x}{c} \, ,x) \right|^2 + \left| \partial_t \, \delta^e_-(t + \frac{x}{c} \, ,x) \right|^2
+
\left|\partial_t \,  \delta^e_+(t - \frac{x}{c} \, ,x)\right|^2 + \left| \partial_t \, \delta^e_-(t - \frac{x}{c} \, ,x) \right|^2
\right]\\   \nonumber
&= 
&\int_{t- \frac{x}{c}}^{t + \frac{x}{c}}
-\frac{1}{c} \, ( \partial_t \, \delta^e_+ + \partial_t \, \delta^e_- )  \,
\partial_t\left(
\sigma^e (\delta_+^e+S_+^e,\,\delta_-^e + S_-^e)
-
\sigma^e (S_+^e,\,  S_-^e)
\right)
\,
\, ds
\\   \nonumber
&  &+  \frac{1}{ c}\left[ \left| \partial_t \, \delta^e_-(t + \frac{x}{c} \, ,x) \right|^2
+
\left|\partial_t \,  \delta^e_+(t - \frac{x}{c} \, ,x)\right|^2 
\right]\\
 \nonumber
&\geq & 
-\int_{t- \frac{x}{c}}^{t + \frac{x}{c}}
\frac{1}{c} \, ( \partial_t \, \delta^e_+ + \partial_t \, \delta^e_- )  \,
\partial_t\left(
\sigma^e (\delta_+^e+S_+^e,\,\delta_-^e + S_-^e)
-
\sigma^e (S_+^e,\,  S_-^e)
\right)
\,
\, ds\, .
\end{eqnarray}
Using the specific form of $\sigma^e$ 
from
(\ref{sigmadefinition}),
the reverse triangle inequality,
 assumption (\ref{druckvor})
and
 (\ref{janbedingung}) 
we obtain 
\begin{eqnarray*}
& &
( \partial_t \, \delta^e_+ + \partial_t \, \delta^e_- )  \,
\partial_t\left(
\sigma^e (\delta_+^e+S_+^e,\,\delta_-^e + S_-^e)
-
\sigma^e (S_+^e,\,  S_-^e)
\right)
\\
& = &
2\, \nu^e ( \partial_t \, \delta^e_+ + \partial_t \, \delta^e_- )  \,
[ | \delta_+^e + S_+^e - \delta_-^e - S_-^e  | 
\partial_t  (\delta_+^e + S_+^e - \delta_-^e - S_-^e )
- | S_+^e - S_-^e  | \, \partial_t (S_+^e  -   S_-^e ) ]
\\
& = & 
2\, \nu^e 
( \partial_t \, \delta^e_+ + \partial_t \, \delta^e_- ) 
\left[ \partial_t  (\delta_+^e  - \delta_-^e  )
| \delta_+^e+S_+^e - \delta_-^e - S_-^e  | \right.
\\
&&\quad \qquad \qquad \left. +
\partial_t  ( S_+^e  - S_-^e )
(| \delta_+^e +S_+^e - \delta_-^e - S_-^e  | - | S_+^e - S_-^e  |)  \right]
\\
& \leq & 
2\, \nu^e 
|\partial_t \, \delta^e_+|^2 \, 
\tilde M 
+
2\, \nu^e \,
|\partial_t \, \delta^e_+ + \partial_t \, \delta^e_- |
\,
\tilde B \, 
| \delta_+^e -   \delta_-^e   | 
\\
& \leq & 
2\, \nu^e 
|\partial_t \, \delta^e_+|^2 \, 
\tilde M 
+
 \nu^e \,
 \tilde B \, 
\left[ |\partial_t \, \delta^e_+ + \partial_t \, \delta^e_- |^2
+
| \delta_+^e -   \delta_-^e   |^2 \right]
\\
& \leq & 
2\,  \nu^e \,
 \tilde B \, 
 \left( | \delta_+^e|^2  +   |\delta_-^e|^2 \right)
+
 2\, \nu^e \,
 \left[ \tilde M +  \tilde B \right] \, 
\left[|\partial_t \, \delta^e_+|^2 + 
|\partial_t \, \delta^e_- |^2
 \right].
\end{eqnarray*}
This yields 
\begin{eqnarray*}
\frac{d}{dx}{\cal K}^e(x) 
&\geq &
- \frac{2\, \nu^e}{c}
\int_{t- \frac{x}{c}}^{t + \frac{x}{c}} 
\tilde B \, 
 \left( | \delta_+^e|^2  +   |\delta_-^e|^2 \right)
+
 \left[ \tilde M +  \tilde B \right]
\left[|\partial_t \, \delta^e_+|^2 + 
|\partial_t \, \delta^e_- |^2
 \right]
 \, dx
\\
&\geq &
- \frac{4\, \nu^e}{c}
 \tilde B \, 
{\cal H}^e(x) 
- \frac{4\, \nu^e}{c}
 \left[ \tilde M +  \tilde B \right]
 {\cal K}^e(x) 
.
\end{eqnarray*}

Thus, using \eqref{eq:deltaboundary}, we obtain
\begin{equation}
    \frac{d}{dx}  ({\cal H}^e(x) +  {\cal K}^e(x)  ) \geq 
    - 4\frac{\nu^e \,  (\tilde M + \tilde B) }{c}  
\big( 
{\cal K}^e(x)
+
 { {\cal H}^e(x) }
\big)
\end{equation}
for $x\in [0,\,L^e]$ almost everywhere.
Due to Gronwall's Lemma
for all $x\in [0,\, L^e]$ 
this implies
\begin{equation}
\label{2201hh1}
{\cal H}^e(x)  +  {\cal K}^e(x)  \leq
\exp \left(4\, \frac{\nu^e}{c} 
(\tilde M + \tilde B) \, (L^e  - x) \right) \, 
\left( {\cal H}^e(L^e)+
{\cal K}^e(L^e) \right).
\end{equation}

%
%
%
%
%
This enables us to estimate
$|\partial_t \delta^e_+(t,x)|^2$
\begin{eqnarray}
&=&\nonumber
\Big|  \partial_t \delta^e_+(t+ \frac{L^e-x}{c},L^e)  + \frac {2 \, \nu^e}c \int_x^{L^e} | \delta^e_+ - \delta^e_- + S^e_+ - S^e_-| \, \partial_t (\delta^e_+ - \delta^e_- + S^e_+ - S^e_-)\\
&& \nonumber \qquad \qquad \qquad \qquad \qquad \qquad \qquad \qquad \qquad \qquad \qquad \qquad 
- | S^e_+ - S^e_- |  \, \partial_t ( S^e_+ - S^e_-)  \, ds \Big|^2 \nonumber\\
& \leq &\nonumber
\Big| 
\big|\partial_t \delta^e_+(t+ \frac{L^e-x}{c},L^e) 
\big|
\\ 
&& \qquad +\frac {2 \, \nu^e}c \int_x^{L^e} 
\left|
\left|\delta^e_+ - \delta^e_- + S^e_+ - S^e_- \right| \, \partial_t (\delta^e_+ - \delta^e_- + S^e_+ - S^e_-)
- | S^e_+ - S^e_- |  \, \partial_t ( S^e_+ - S^e_-)
\Big| \,  ds \right|^2 \nonumber\\
&  \leq & \nonumber
\Bigg| 
\big|\partial_t \delta^e_+(t+ \frac{L^e-x}{c},L^e) 
\big|\\
&& \qquad
+ \frac {2 \, \nu^e}c \int_x^{L^e} 
\left|\delta^e_+ - \delta^e_- \right| \, 
\left|\partial_t (S^e_+ - S^e_-)\right|
+ 
\left|
\delta^e_+ - \delta^e_- + S^e_+ - S^e_- \right| 
\, \left|\partial_t ( \delta^e_+ - \delta^e_-) 
\right|\, ds 
\Bigg|^2 \nonumber\\
& \leq &
\nonumber
2\, \big|\partial_t \delta^e_+(t+ \frac{L^e-x}{c},L^e) 
\big|^2\\
&& \nonumber \qquad 
+ \frac{8 \, (\nu^e)^2}{c^2}
\left|
\int_x^{L^e} 
\left|\delta^e_+ - \delta^e_- \right| \, 
\left|\partial_t (S^e_+ - S^e_-)\right|
+ 
\left|
\delta^e_+ - \delta^e_- + S^e_+ - S^e_- \right| 
\, \left|\partial_t ( \delta^e_+ - \delta^e_-) 
\right|\, ds 
\right|^2
\\
&\leq &
\nonumber
2\, 
\Big|  \partial_t \delta^e_+(t+ \frac{L^e-x}{c},L^e) \Big|^2 
 + 
  \frac{8 \, (\nu^e)^2}{c^2}
\left|
\int_x^{L^e} 
\tilde B
\left|\delta^e_+ - \delta^e_- \right| \, 
+ 
\tilde M
\, \left|\partial_t ( \delta^e_+ - \delta^e_-) 
\right|\, ds 
\right|^2
\\
&\leq &
\nonumber
2\, 
\Big|  
\partial_t \delta^e_+(t+ \frac{L^e-x}{c},L^e) \Big|^2 
 + 
  \frac{16 \, (\nu^e)^2\, L^e}{c^2}
\,
\int_x^{L^e} 
(\tilde B)^2
\left|\delta^e_+ - \delta^e_- \right|^2 \, 
+ 
(\tilde M)^2
\left|
\partial_t ( \delta^e_+ - \delta^e_-) 
\right|^2\, ds 
\end{eqnarray}
where the integrands are to be evaluated at $(t + \tfrac{s-x}{c},s).$
Thus we have

\begin{eqnarray}
&& \nonumber \int_0^{L^e}|\partial_t \delta^e_+(t,x)|^2 \, dx\\
\nonumber 
& \leq &
2\, 
\int_0^{L^e}
\Big|  
\partial_t \delta^e_+(t+ \frac{L^e-x}{c},L^e) \Big|^2 
\, dx
\\
\nonumber
&  +  &
   \frac{16 \, (\nu^e)^2\, L^e}{c^2}
\int_0^{L^e}
\,
\int_x^{L^e} 
(\tilde B)^2
\left|\delta^e_+ - \delta^e_- \right|^2 \, 
+ 
(\tilde M)^2
\left|
\partial_t ( \delta^e_+ - \delta^e_-) (t + \tfrac{s-x}{c},s)
\right|^2\, ds \, dx
\\
\nonumber
& \leq &
2\, 
\int_0^{L^e}
\Big|  
\partial_t \delta^e_+(t+ \frac{L^e-x}{c},L^e) \Big|^2 
\, dx
\\
\nonumber
&  +  &
   \frac{16 \, (\nu^e)^2\, L^e}{c}
\int_0^{L^e}
\,
\int_{t - \frac{x}{c}}
^{t + \frac{x}{c}} 
(\tilde B)^2
\left|(\delta^e_+ - \delta^e_-)
(s,x)\, \right|^2 
+ 
(\tilde M)^2
\left|
\partial_t ( \delta^e_+ - \delta^e_-) 
(s,x)
\right|^2\, ds \, dx.
\end{eqnarray}
Hence   using definition (\ref{hdefinition}) of  ${\cal H}^e(x)$
and
definition (\ref{kdefinition}) of  ${\cal K}^e(x)$
we obtain
\begin{multline}
\int_0^{L^e}|\partial_t \delta^e_+(t,x)|^2 \, dx
 \leq 
2\, 
\int_0^{L^e}
\Big|  
\partial_t \delta^e_+(t+ \frac{L^e-x}{c},L^e) \Big|^2 
\, dx
\\
  +  
   \frac{64 \, (\nu^e)^2\, L^e}{c}
\int_0^{L^e}
\,
(\tilde B)^2
{\cal H}^e(x)
+ 
(\tilde M)^2
{\cal K}^e(x)
 \, dx\, .
\end{multline}
With  (\ref{2201hh1})
and the notation 
$\Upsilon =  \frac{\max\{(\tilde B)^2,\, (\tilde M)^2\} }
   { \tilde M + \tilde B  }
$
this yields
\begin{eqnarray*}
&&\int_0^{L^e}|\partial_t \delta^e_+(t,x)|^2 \, dx\\
\nonumber 
& \leq &
2\, 
\int_0^{L^e}
\Big|  
\partial_t \delta^e_+(t+ \frac{L^e-x}{c},L^e) \Big|^2 
\, dx
\\
\nonumber
&  +  &
   \frac{64 \, (\nu^e)^2\, L^e}{c}
   \max\{(\tilde B)^2,\, (\tilde M)^2\}
\int_0^{L^e}
\exp \left(4\, \frac{\nu^e}{c} 
(\tilde M + \tilde B) \, (L^e  - x) \right) \, 
\left( {\cal H}^e(L^e)+
{\cal K}^e(L^e) \right)
\, dx
\\
& \leq &
2\, 
\int_0^{L^e}
\Big|  
\partial_t \delta^e_+(t+ \frac{L^e-x}{c},L^e) \Big|^2 
\, dx
\\
\nonumber
&  +  &
   \left(16 \, \nu^e\, L^e\right)
  \Upsilon
\exp \left(4\, \frac{\nu^e}{c} 
(\tilde M + \tilde B) \, L^e  \right) \, 
\left( {\cal H}^e(L^e)+
{\cal K}^e(L^e) \right)
\\
& \leq &
2 \, c\, 
\int_t^{ t + \frac{L^e}{c}}
\left| \partial_t \delta^e_+(s, \, L^e) \right|^2\, ds
\\
&
  +
&  
    \left(16 \, \nu^e\, L^e\right)
  \Upsilon
\exp \left(4\, \frac{\nu^e}{c} 
(\tilde M + \tilde B) \, L^e  \right) \, 
\left( {\cal H}^e(L^e)+
{\cal K}^e(L^e) \right).
\end{eqnarray*}
Similarly we obtain
the same upper bound for
$
\int_0^{L^e}|\partial_t \delta^e_-(t,x)|^2 \, dx.
$

Adding up the inequalities for $\partial_t\delta^e_+$ and $\partial_t\delta^e_-$ yields
\begin{eqnarray*}
 & & \int_0^{L^e}|\partial_t \delta^e_+(t,x)|^2
+
|\partial_t \delta^e_-(t,x)|^2
\, dx
\\
& \leq &
4 \, c\, {\cal K}^e(x^e(v))
  +
    \left(32 \, \nu^e\, L^e\right)
  \Upsilon \,
\exp \left(4\, \frac{\nu^e}{c} 
(\tilde M + \tilde B) \, L^e  \right) \, 
\left( {\cal H}^e(x^e(v))+
{\cal K}^e(x^e(v)) \right).
\end{eqnarray*}
For the case $x^e(v)=0$ we obtain 
the analogous inequality.
Adding up this inequality and (\ref{observabilityinequality})
yields the desired observability inequality
(\ref{observabilityinequalityh1}) with
\[
C_1(\tilde M,\, \tilde B) 
= C_0(\tilde M ) + 
2 \, c
+
\max_{e\in E}
 \left(16 \, \nu^e\, L^e\right)
\Upsilon \, 
\exp \left(4\, \frac{\nu^e}{c} 
(\tilde M + \tilde B) \, L^e  \right) .
\]




\section{Exponential decay of the observer error  on the network}
\label{Sec:Exponential_decay}

In this section, we  analyse the
evolution of the state $R$ of
the observer   ${\bf (R)}$.
We show that  $R$ 
approaches the  state $S$  of system ${\bf (S)}$  exponentially fast on each edge.
In order to show this, we study the evolution of
the error system ${\bf (Diff)}$ and
show that the solution $\delta$
decays exponentially fast. 

Theorem \ref{stability} has two parts.
In the first part, a sufficient condition 
for the exponential decay of the $L^2$-norm  (\ref{shiri}) 
on a finite time interval $[0, \, T]$ 
is provided under the  assumption (\ref{druckvor}).
This first part of Theorem \ref{stability} can 
be applied to $L^\infty$-solutions
as discussed in Theorem \ref{existence}.
For the proof of the first part, 
the observability inequality from Theorem \ref{observabilitythm} is used.

In the second part of Theorem \ref{stability}, more
regular $H^1$-solutions are considered.
For the proof,  
the observability inequality from Theorem \ref{observabilitythmh1} is used.

\begin{theorem}
\label{stability}
Define $T_0 = \max_{e\in E}  \frac{L^e}{c}$.
Let $T > 2 \, T_0$  be given.
For all $e\in E$, let initial states 
$y^e_+$,  $y^e_-$, $z^e_+$, $z^e_-
\in  L^\infty(0, \, L^e)$ be given.
Assume that for each  node
$v\in V$ a number 
$|\mu^v|  \leq 1$
is given.

Assume that there exists a set $ \tilde V \subset V$ with the following property:
For all $e\in E$ there exists
$v\in \tilde V$ such that
$e\in E_0(v)$ and 
$ |\mu^v|<1 $.

Assume that 
there exists a real number $J$ such that 
for the initial states 
$y^e_\pm - J$, $z^e_\pm - J$
and for  the  boundary  controls
(for all $e\in E$ where $|E_0(v)|=1$) 
$
u^e -J \in L^\infty(0, \, T)
$
have  a sufficiently small $L^\infty$-norm
such that  solutions of  systems ${\bf (S)}$
and 
${\bf (R)}$
exist on $[0, \, T]$
in $\times_{e\in E}   L^\infty((0,\, T) \times (0, \, L^e))$
and satisfy (\ref{druckvor}). 
Then  the  solution of system ${\bf (Diff)}$
%
is   exponentially stable in the sense that
there exist  constants
$C_1 >0$ and $\mu_0>0$
such that for all $t\in [0,\, T]$ 
the following inequality holds: 
\begin{multline}
\label{shiri}
\sum_{e\in E}
 (D^e)^2 \int_0^{L^e} \left|\delta^e_+(t,\, x)  \right|^2 
 + \left| \delta^e_-(t,\, x) \right|^2
\, dx\\
\leq C_1
\exp(-\mu_0 \,  t) \,
\sum_{e\in E}
 (D^e)^2 \int_0^{L^e} \left|\delta^e_+(0,\, x)  \right|^2
 + \left| \delta^e_-(0,\, x) \right|^2 \, dx.
\end{multline}
Hence the $L^2$-norm of the difference $\delta$ between the  state $R$
of the observer and
the  state  $S$ of the original system  decays exponentially fast.

 Assume in addition that
$\partial_t y^e_\pm$
has a sufficiently small $L^\infty$-norm
and is compatible with the node condition and the boundary conditions
such that a solution of  system ${\bf (S)}$ exists on $[0, \, T]$
in $\times_{e\in E}  C([0,\, T],\, W^{1,\infty}(0, \, L^e))$
and satisfies (\ref{druckvor})
and
(\ref{janbedingung}).
Assume that
\begin{equation}
    \label{observiassumption}
   8  \, T_0\,  \left(\Delta(\nu T_0)\right)^2 \,  \max_{e\in E} \nu^e 
\leq 
\frac{c }{C_1( \tilde M, \, \tilde B)}  
\min_{e \in E} \left[
\sum_{v \in V_0(e)}  \frac{1-|\mu^v|^2}{1+|\mu^v|^2}
\right]
\end{equation} 
where
\begin{equation}
    \label{observiassumption2}\Delta(\nu T_0) := 
  \exp\left( 8 \,\max_{e\in E} \nu^e
    \tilde B   \, T_0 \right)\,.\end{equation}

   Then 
   in addition to (\ref{shiri})
  also  the $L^2$-norm of the  time-derivatives decay exponentially fast in the sense that
  there exists a number
  $\mu_1 \in (0,\, \mu_0)$
  such that for some $\tilde C>0$ 
  we have 
    \begin{multline*}
   \sum_{e\in E}
 (D^e)^2   \int_0^{L^e}\left|  \partial_t \delta^e_+(t,\, x) \right|^2 +
\left|  \partial_t  \delta^e_-(t,\, x) \right|^2
 \, dx\\
 \leq
 \frac{\tilde C}{
 \mbox{\textrm{e}}^{\mu_1  \, t }}
  \;\sum_{e\in E}
  (D^e)^2   
  \Bigg[
  \int_0^{L^e}\left|  \partial_t \delta^e_+(0,\, x) \right|^2 +
\left|  \partial_t  \delta^e_-(0,\, x) \right|^2
+ 
\left|\delta^e_+(0,\, x)  \right|^2
 + \left| \delta^e_-(0,\, x) \right|^2
 \, dx
\Bigg].
  \end{multline*}

\end{theorem}

\begin{remark}
It is important to note that for the exponential decay of the $L^2$-norm we do not
need any restrictions on the lengths of $L^e$.

Assumption \eqref{observiassumption}   holds if
for all $ e\in E$, friction coefficients $\nu^e$ 
or  the lengths $L^e$ are sufficiently small.
This is similar to the assumptions in \cite{guherty}
where the decay of the $L^2$-norm has been studied.
%
Note that assumption
 \eqref{observiassumption}  
also holds
if the initial data have a sufficiently small $H^1$-norm
such that we can choose $\tilde M$ and $\tilde B$ sufficiently small.
This implies that  for arbitrary large values of the lengths $L^e$,
there exists an  $H^1$-neighbourhood of
the constant state, 
such for initial states in this neighbourhood 
 the $H^1$-norm of  the system state decays exponentially.

\end{remark}

\begin{remark}\label{janremark_decay}
We can prove a similar result to Theorem \ref{stability}
 if we replace assumption \eqref{janbedingung} by \eqref{janbesser}. 
 Going to the weaker assumption \eqref{janbesser} reduces the decay rates $\mu_0$ and $\mu_1$ and leads to a much more restrictive condition on the problem parameters in place of \eqref{observiassumption}.
 \end{remark}

{\bf Proof of Theorem \ref{stability}.}
%
%
%
%
Let $\bar t\in (0, t) $ be given. 
For $e\in E$ the partial differential equation in  ${\bf (Diff)}$ implies that 
\[\partial_t  \delta^e_+(t,\, x)  
=
-c \, \partial_x  \delta^e_+(t,\, x)  - 
\left[\sigma^e(\delta^e_+ + S^e_+,\, \delta^e_- + S^e_-)
-
\sigma^e(S^e_+,\, S^e_-)
\right].
\]
We multiply this equation by $\delta^e_+ $  and integrate
over the interval
$(t - \bar t,\, t + \bar t)$  to obtain
\begin{multline}
 \int_{ t - \bar t}^{t + \bar t  }
 \int_0^{L^e}
 \delta^e_+(\tau,\, x)  \,
 \partial_t  \delta^e_+(\tau,\, x) \, dx \, d\tau\\
 =
 \int_{ t - \bar t}^{t +  \bar t}
 \int_0^{L^e} 
 -c \, \delta^e_+(\tau,\, x)   \, \partial_x  \delta^e_+(t,\, x)    \\
 - \delta^e_+(\tau,\, x)   
  \left[\sigma^e(\delta^e_+ + S^e_+,\, \delta^e_- + S^e_-)
-
\sigma^e(S^e_+,\, S^e_-)
\right]
 \, dx \, d\tau.
\end{multline}
This yields
\begin{multline}
- \frac{c}{2}  \, \int_{ t - \bar t}^{t + \bar t} \left[  |\delta^e_+(\tau,\, x)|^2 \right]|_{x=0}^{L^e}\\
-  
 \int_0^{L^e} 
   \delta^e_+(\tau,\, x) \, 
   \left[\sigma^e(\delta^e_+ + S^e_+,\, \delta^e_- + S^e_-)
-
\sigma^e(S^e_+,\, S^e_-)
\right]
   \, dx \, d\tau\\
=
\frac{1}{2} \,  \int_0^{L^e} \left[|\delta^e_+(\tau,\, x)|^2  \right]|_{\tau  = t-  \bar t}^{t+  \bar t}\, dx .
\end{multline}
Similarly, we obtain
\begin{multline}
 \frac{c}{2}  \, \int_{ t - \bar t}^{t +\bar t } \left[  |\delta^e_-(\tau,\, x)|^2 \right]|_{x=0}^{L^e}
\\+
 \int_0^{L^e} 
   \delta^e_-(\tau,\, x)  
   \left[\sigma^e(\delta^e_+ + S^e_+,\, \delta^e_- + S^e_-)
-
\sigma^e(S^e_+,\, S^e_-)
\right]
   \, dx \, d\tau\\
=
\frac{1}{2} \,  \int_0^{L^e} \left[|\delta^e_-(\tau,\, x) |^2  \right]|_{\tau  = t-  \bar t }^{t+  \bar t}\, dx .
\end{multline}
For $e\in E$ and $t \in [0, T]$, we define  
\begin{equation}
\label{Ledefinition}
{\cal L}_0^e(t) := \frac{(D^e)^2}{2} \,  \int_0^{L^e} \left|\delta^e_+(t,\, x)  \right|^2 + \left|\delta^e_-(t,\, x)
\right|^2
 \, dx
\ \text{and} \
{\cal L}_0(t) :=  \sum_{e\in E} {\cal L}_0^e(t).
\end{equation}
Then we have
\begin{multline*}{\cal L}_0(t  + \bar t) -   {\cal L}_0(t - \bar t ) 
=
\sum_{e\in E}
- \frac{c \, (D^e)^2}{2}  \, \int_{ t - \bar t}^{t +  \bar t} \left[  |\delta^e_+(\tau,\, x)  |^2 
-   |\delta^e_-(\tau,\, x)  |^2 \right]\Big|_{x=0}^{L^e}
\\
-  
 \int_0^{L^e}  (D^e)^2
   (\delta^e_+(\tau,\, x)  - \delta^e_-(\tau,\, x) ) \,
   \left[\sigma^e(\delta^e_+ + S^e_+,\, \delta^e_- + S^e_-)
-
\sigma^e(S^e_+,\, S^e_-)
\right] \,
   dx \, d\tau.
\end{multline*}
Hence due to the definition of $\sigma^e$ 
with the mean value theorem we obtain the inequality
\[{\cal L}_0(t  + \bar t ) -   {\cal L}_0(t -  \bar t   ) 
\leq
\sum_{e\in E}
- \frac{c \, (D^e)^2}{2}  \, \int_{ t - \bar t}^{t +  \bar t} \left[  |\delta^e_+(\tau,\, x) |^2 
-   |\delta^e_-(\tau,\, x)|^2 \right]\Big|_{x=0}^{L^e} \, d\tau.
\]
i.e.
\begin{multline}{\cal L}_0(t  + \bar t ) -   {\cal L}_0(t -  \bar t   ) 
\\
\leq
\sum_{e\in E}
 \frac{c \, (D^e)^2}{2}  \,
 \int_{ t - \bar t}^{t +  \bar t}  
 -
  |\delta^e_+(\tau,\, L^e) |^2 
+   |\delta^e_-(\tau,\, L^e)|^2 
 +
 |\delta^e_+(\tau,\, 0) |^2 
-   |\delta^e_-(\tau,\, 0)|^2  \, d\tau.
\end{multline}

At any interior node $v\in V$   (i.e. $|E_0(v)|>1$)
the node conditions  imply 
\begin{equation}
    \label{middlenode_pl3537}
\sum_{e\in E_0(v)}
  (D^e)^2 \, |\delta^e_{out}(\tau,\, x^e(v)) |^2 
= |\mu^v|^2 \sum_{e\in E_0(v)}
  (D^e)^2 \,   |\delta^e_{in}(\tau,\, x^e(v))|^2, 
  \end{equation}
 hence  we have 
\begin{multline}
    \label{middlenode}
  \sum_{e\in E_0(v)}  (D^e)^2 \,
 \left[  |\delta^e_{out}(\tau,\, x^e(v)) |^2 -  |\delta^e_{in}(\tau,\, x^e(v)) |^2\right]
\\
\leq \frac{ |\mu^v|^2 - 1}{|\mu^v|^2+1}
\sum_{e\in E_0(v)}
  (D^e)^2 \,  \left[
  |\delta^e_-(\tau,\, x^e(v))|^2  +  |\delta^e_+(\tau,\, x^e(v))|^2\right].
\end{multline}
Similarly, the boundary conditions at 
any boundary node 
$v\in V$ (with $|E_0(v)|=1$)
imply
\begin{equation}
    \label{boundarynode_pl}
   |\delta^e_{out}(\tau,\, x^e(v)) |^2 
= |\mu^v|^2   
|\delta^e_{in}(\tau,\, x^e(v))|^2,  \text{ for } e \in E_0(v)
  \end{equation}
  so that
\begin{multline}
    \label{boundarynode}
  \sum_{e\in E_0(v)}  (D^e)^2 \,
 \left[ - |\delta^e_{in}(\tau,\, x^e(v)) |^2 +  |\delta^e_{out}(\tau,\, x^e(v)) |^2\right]
\\
\leq \frac{ |\mu^v|^2 - 1}{|\mu^v|^2+1}
\sum_{e\in E_0(v)}
  (D^e)^2 \,  \left[
  |\delta^e_-(\tau,\, x^e(v))|^2  +  |\delta^e_+(\tau,\, x^e(v))|^2\right].
\end{multline}
This yields
\begin{multline}\label{eq:interm}
{\cal L}_0(t  + \bar t ) -   {\cal L}_0(t -  \bar t   ) \\
\leq
\sum_{v \in V}
\frac{ |\mu^v|^2 - 1}{|\mu^v|^2+1}
 \sum_{e\in E_0(v)} 
\frac{c \, (D^e)^2}{2}
\int_{ t - \bar t}^{t +  \bar t} 
 |\delta^e_-(\tau,\, x^e(v)) |^2 +
 |\delta^e_+(\tau,\, x^e(v)) |^2
 \, d\tau\,.
\end{multline}

Since $|\mu^v|^2 \leq 1$ for all $v \in V$  ,
in particular, we have
${\cal L}_0(t  + \bar t) \leq {\cal L}_0(t - \bar t  )$.
Since the above inequality can be derived for all $\bar t \in (0, t)$,
 this implies 
that  ${\cal L}_0$ is decreasing. 

We choose $\bar t =  T_0>0$.
For all $v\in V$ and $e \in E_0(v)$ the observability inequality  \eqref{observabilityinequality}  implies
\begin{multline}
\int_{t - T_0 }^{t + T_0} \left| \delta^e_+(\tau,\, x^e(v)) \right|^2 +  \left| \delta^e_-(\tau,\, x^e(v))  \right|^2 \, d\tau \\
\geq \frac{1}{ C_0(\tilde M)} \,
 \int_0^{L^e}
|\delta^e_+(t,\, x)|^2 + |\delta^e_-(t,\, x)|^2 \, dx
\end{multline}
with $C_0(\tilde M)$ as defined in (\ref{c0definition}).
Inserting this into \eqref{eq:interm}
implies

\begin{multline}
    {\cal L}_0(t + \bar t) -   {\cal L}_0(t - \bar t) \\
    \leq
    \frac{1}{C_0(\tilde M)}\sum_{e\in E} \sum_{v \in V_0(e)}  \frac{|\mu^v|^2- 1}{|\mu^v|^2+ 1}   \frac{c (D^e)^2}{2}   \int_0^{L^e}
|\delta^e_+(t,\, x)|^2 + |\delta^e_-(t,\, x)|^2 \, dx
\end{multline}
where $V_0(e) $ denotes the set of nodes adjacent to $e$.
This yields the inequality
\[{\cal L}_0(t  + T_0 ) -   {\cal L}_0(t - T_0 )
\leq
-
\frac{c}{  C_0(\tilde M) }  
\min_{e\in E}
\left\{ 
\sum_{v \in V_0(e)}
\frac{1- |\mu^v|^2 }{|\mu^v|^2+1}
\right\}
 {\cal L}_0(t ).
\]
Define the constant
\[\upsilon_0:=
\min_{e\in E}
\left\{ 
\sum_{v \in V_0(e)}
\frac{1- |\mu^v|^2 }{|\mu^v|^2+1}
\right\}
\]
Since ${\cal L}_0$ is decreasing this yields
\[{\cal L}_0(t  + T_0 )
\leq
 {\cal L}_0(t - T_0  )
-
\frac{c}{ C_0(\tilde M)  } \, 
\upsilon_0
{\cal L}_0(t + T_0 ).
\]
Hence we have
\[{\cal L}_0(t  + T_0 )
\leq
\frac{1}{1 + \frac{c}{ C_0(\tilde M)}  \, 
\upsilon_0
} \, {\cal L}_0(t - T_0  ).
\]
Similarly as in  Lemma 2 from \cite{gugattucsnak}, this implies that ${\cal L}_0$ decays exponentially fast.

Now we consider the evolution of the time-derivatives.
Note that here the analysis is more involved than
for the $L^2$-estimate, since the sign
of the friction term cannot be determined a priori.
In addition, the 
time derivative of the friction term
requires an estimate for the term
$\partial_t (S^e_+ - S^e_-)$,
see assumption (\ref{janbedingung}).

Similar as in the proof of the observability inequality,
it is necessary to consider the sum
${\cal L}_0 + {\cal L}_1$
to obtain an estimate.

For $e\in E$ and $t\in [0, T]$, we define
\begin{equation}
\label{Ledefinitionobserv}
{\cal L}_1^e(t) := \frac{1}{2} \, (D^e)^2 \int_0^{L^e} \left|  \partial_t \delta^e_+(t,\, x) \right|^2 +
\left|  \partial_t  \delta^e_-(t,\, x) \right|^2
 \, dx \ \text{and} \
{\cal L}_1(t) :=  \sum_{e\in E} {\cal L}_1^e(t).
\end{equation}
Due to the partial differential equation in system ${\bf (S)}$,
for solutions with $H^2$-regularity we have
\[
\partial_{tt} \left(
\begin{array}{r}
\delta_+^e
\\
\delta_-^e
\end{array}
\right)
=
c \, \partial_{xt} 
\left(
\begin{array}{r}
-  \delta_+^e
\\
  \delta_-^e
\end{array}
\right)
+
\partial_{t}
\left(  
\sigma^e(\delta^e_+ + S^e_+, \,
\delta^e_- + S^e_-) -
\sigma^e( S^e_+, \, S^e_-)
\right)
\;
\left(
\begin{array}{r}
- 1
\\
1
 \end{array}
\right).
\]
For the time-derivative of ${\cal L}_1^e$ we have
\begin{eqnarray*}
\frac{d}{dt}  {\cal L}^e_1(t) & = &
 (D^e)^2 \int_0^{L^e}  \, \partial_t \delta^e_+(t,\, x) \, \partial_{tt} \delta^e_+(t,\, x)
    +   \partial_t  \delta^e_-(t,\, x)  \, \partial_{t t}  \delta^e_-(t,\, x)  \, dx
    .
\end{eqnarray*}
With the partial differential equation for $\partial_t \delta$
this yields 
\begin{eqnarray*}
\frac{d}{dt}  {\cal L}^e_1(t) 
 & = &
 (D^e)^2 \int_0^{L^e}
   \partial_t \delta^e_+  \, [- c \,\partial_{tx} \delta^e_+
-
\partial_t \left(  
\sigma^e(\delta^e_+ + S^e_+, \,
\delta^e_- + S^e_-) -
\sigma^e( S^e_+, \, S^e_-)
\right)
]
 \\
 &
  &
\qquad \qquad+\partial_t \delta^e_- \, [c \, \partial_{tx} \delta^e_- + 
 \partial_t \left(  
\sigma^e(\delta^e_+ + S^e_+, \,
\delta^e_- + S^e_-) -
\sigma^e( S^e_+, \, S^e_-)
\right)]
 \, dx 
 \\
  & = &
 (D^e)^2 \int_0^{L^e}
   \partial_t  \delta^e_+  \, [-  c \,\partial_{xt}\delta^e_+  ]
    +   \partial_t \delta^e_-(t,\, x) \, [ c \, \partial_{tx} \delta^e_-  ]
 \\
 & &
   \qquad\qquad-   
   \left[
   \partial_t \delta_+^e -
  \partial_t \delta_-^e
  \right]\, 
  \partial_{ t} \left(  
\sigma^e(\delta^e_+ + S^e_+, \,
\delta^e_- + S^e_-) -
\sigma^e( S^e_+, \, S^e_-)
\right)
   \, dx
    .
\end{eqnarray*}

This yields
\begin{eqnarray*}
\frac{d}{dt}  {\cal L}^e_1(t) 
 & =  &
c \, \frac{ (D^e)^2 }{2} \int_0^{L^e} -    \partial_x \left( \partial_t \delta^e_+(t,\, x)  \right)^2
+
   \partial_x \left( \partial_t \delta^e_-(t,\, x)  \right)^2
   \, dx
  \\
 &  &
 -
   (D^e)^2 \int_0^{L^e}
 \partial_t  \left(  
\sigma^e(\delta^e_+ + S^e_+, \,
\delta^e_- + S^e_-) -
\sigma^e( S^e_+, \, S^e_-)
\right) \,  \left[
  \partial_t \delta_+^e -
   \partial_t \delta_-^e
    \right]\, dx.
\end{eqnarray*}
Note that we have
\begin{eqnarray*}
& &
\partial_t \left(  
\sigma^e(\delta^e_+ + S^e_+, \,
\delta^e_- + S^e_-) -
\sigma^e( S^e_+, \, S^e_-)
\right) \,  
\left[
  \partial_t \delta_+^e -
   \partial_t \delta_-^e
    \right]
\\
&
=
&
2\,\nu^e \, 
 \left[ 
 | \delta^e_+ - \delta^e_- + S^e_+ - S^e_-| \, \partial_t (\delta^e_+ - \delta^e_- + S^e_+ - S^e_-)
- | S^e_+ - S^e_- |  \, \partial_t ( S^e_+ - S^e_-)
\right]
\left[
  \partial_t \delta_+^e -
   \partial_t \delta_-^e
    \right]
    \\
       & = &
    2\,\nu^e \, 
 \left\{
\left[ 
 | \delta^e_+ - \delta^e_- + S^e_+ - S^e_-| 
 -
 | S^e_+ - S^e_- |
 \right]
 \partial_t ( S^e_+ - S^e_-)\right.
\\
&& \qquad \quad \qquad \quad \qquad \quad \qquad +\left.
 | \delta^e_+ - \delta^e_- + S^e_+ - S^e_-|
 \,
  \partial_t ( \delta^e_+ - \delta^e_-)
\right\}
\left[
  \partial_t \delta_+^e -
   \partial_t \delta_-^e
    \right].
\end{eqnarray*}
Due to
(\ref{druckvor})
and
(\ref{janbedingung})
this yields the lower bound
\begin{eqnarray*}
& &
\Big[ \partial_t \left(  
\sigma^e(\delta^e_+ + S^e_+, \,
\delta^e_- + S^e_-) -
\sigma^e( S^e_+, \, S^e_-)
\right) \,  
\left[
  \partial_t \delta_+^e -
   \partial_t \delta_-^e
    \right]
    \Big]
\\
&
\geq 
&
2\,\nu^e \, 
 | \delta^e_+ - \delta^e_- + S^e_+ - S^e_-| \,
\left[
  \partial_t \delta_+^e -
   \partial_t \delta_-^e
    \right]^2
    \\
    & &
    -
 2\,\nu^e \, 
 \left\{
 | \delta^e_+ - \delta^e_-  |
 \,
 |
 \partial_t ( S^e_+ - S^e_-)
 |
\,
\left|
  \partial_t \delta_+^e -
   \partial_t \delta_-^e
    \right|
    \right\}
 \\
&
\geq
&
- \nu^e \, 
 \left[
 | \delta^e_+ - \delta^e_-  |^2
 +
 \left|
  \partial_t \delta_+^e -
   \partial_t \delta_-^e
    \right|^2
 \right]
 |
 \partial_t ( S^e_+ - S^e_-)
 |
 \\
&
\geq
&
- \nu^e \, 
 \tilde B
 \,
 \left[
 | \delta^e_+ - \delta^e_-  |^2
 +
 \left|
  \partial_t \delta_+^e -
   \partial_t \delta_-^e
    \right|^2
 \right]
.
\end{eqnarray*}

Integration  yields
\begin{eqnarray*}
 \frac{d}{dt}{\cal L}^e_1(t)
& = &
   c \,  \frac{ (D^e)^2}{2}  \,  \left[  - \left( \partial_t \delta^e_+(t,\, x)  \right)^2
   +
   \left(  \partial_t  \delta^e_-(t,\, x) \right)^2 \right] |_{x=0}^{L^e}
  \\
 & &
  - 
   (D^e)^2 \int_0^{L^e}
  \left(  
\sigma^e(\delta^e_+ + S^e_+, \,
\delta^e_- + S^e_-) -
\sigma^e( S^e_+, \, S^e_-)
\right)_t \,  \left[
   \partial_t \delta_+^e -
   \partial_t \delta_-^e
    \right]\,
   \, dx
     \\
     &  \leq &
 c \,  \frac{ (D^e)^2}{2}  \,  \left[   -   \left( \partial_t \delta^e_+(t,\, L^e)  \right)^2 +     \left(  \partial_t  \delta^e_-(t,\, L^e) \right)^2
                   + \left( \partial_t \delta^e_+(t,\, 0)  \right)^2
            - \left(  \partial_t  \delta^e_-(t,\, 0) \right)^2                  \right]
 \\
&
  &
  +(D^e)^2 \,\nu^e
  \int_0^{L^e}
 \tilde B
 \left[
 | \delta^e_+ - \delta^e_-  |^2
 +
 \left|
  \partial_t \delta_+^e -
   \partial_t \delta_-^e
    \right|^2
 \right]
\, dx
  \\
  &
  \leq
  &
 c \,  \frac{ (D^e)^2}{2}  \,  \left[ 
  -   \left( \partial_t \delta^e_+(t,\, L^e)  \right)^2 +    \left(  \partial_t  \delta^e_-(t,\, L^e) \right)^2
  +
                   \left( \partial_t \delta^e_+(t,\, 0)  \right)^2
            - \left(  \partial_t  \delta^e_-(t,\, 0) \right)^2                  \right]
 \\
&
  &
  +
   4 \,\nu^e
   \tilde B  \,
   \left[
   {\cal L}^e_0(t)
   +
   {\cal L}^e_1(t)
   \right].
\end{eqnarray*}
The boundary conditions and 
 the coupling conditions imply that at any node $v \in V$ we have
\[
\sum_{e\in E} (D^e)^2 (\partial_t \delta^e_{out}(t,x^e(v)))^2 =  
(\mu^v)^2 \sum_{e\in E} (D^e)^2 \, (\partial_t \delta^e_{in}(t,x^e(v)))^2 
\]

 Thus, we obtain the inequality
 \begin{eqnarray}\label{H2H1}
\frac{d}{dt}  {\cal L}_1(t)
& \leq &\nonumber
- \sum_{e\in E}
 c \,  \frac{ (D^e)^2}{2}  \,  \left[ 
 \sum_{v \in V_0(e)}
    \frac{1-|\mu^v|^2}{1 + |\mu^v|^2}\left( \left( \partial_t \delta^e_+(t,\, x^e(v))  \right)^2 +     \left(  \partial_t  \delta^e_-(t,\, x^e(v)) \right)^2\right)
             \right]
 \\
&
  &
 \qquad \quad +   4 \,\nu^e
   \tilde B  \,
   \left[
   {\cal L}^e_0(t)
   +
   {\cal L}^e_1(t)
   \right]
  .
   \end{eqnarray}
   Note that there is no second order derivative in \eqref{H2H1} and it can be extended to $H^1$ solutions by  a density argument.
   As in the $H^1$-norm observability result, we need to consider the sum ${\cal L}_\Sigma:= {\cal L}_0 + {\cal L}_1$ and obtain
   \begin{eqnarray*}
&&\frac{d}{dt} {\cal L}_\Sigma(t)\\
& \leq &
- \sum_{e\in E}
 c \,  \frac{ (D^e)^2}{2}  \,  \Bigg[ 
  \sum_{v \in V_0(e)}  \frac{1-|\mu^v|^2}{1 + |\mu^v|^2}\Bigg( 
    \left( \partial_t \delta^e_+(t,\, x^e(v))  \right)^2 +     \left(  \partial_t  \delta^e_-(t,\, x^e(v)) \right)^2\\
    && \qquad \qquad \qquad \qquad
    + \left(\delta^e_+(t,\, x^e(v))  \right)^2 +     \left(   \delta^e_-(t,\, x^e(v)) \right)^2
    \Bigg)
             \Bigg]
             +  4 \,\nu^e
   \tilde B  \,
   \left[
   {\cal L}^e_0(t)
   +
   {\cal L}^e_1(t)
   \right]
  .
   \end{eqnarray*}

   We integrate over $[t-T_0,t+T_0]$ and obtain
   \begin{eqnarray*}
& &  {\cal L}_\Sigma(t+ T_0) -  {\cal L}_\Sigma(t - T_0)
 \\
& \leq &
 -
   \sum_{v \in V}
     \frac{1-|\mu^v|^2}{1 + |\mu^v|^2}
 \sum_{e\in E_0(v)}
 c \,  \frac{ (D^e)^2}{2}  \, 
 \int_{t- T_0}^{t + T_0}
   \left| \delta^e_+(\tau,\,x^e(v))  \right|^2 +    
   \left|    \delta^e_-(\tau,\, x^e(v)) \right|^2\\
   && \qquad \qquad \qquad \qquad \qquad \qquad \qquad \qquad
  +     \left| \partial_t \delta^e_+(\tau,\, x^e(v))  \right|^2 +     \left|  \partial_t  \delta^e_-(\tau,\, x^e(v)) \right|^2
     \, d\tau
     \\
& &+
  4\,  \max_{e\in E} \nu^e \,\tilde B  \,
  \int_{t- T_0}^{t + T_0}
  {\cal L}_\Sigma(\tau)
   \, d\tau\\
   & \leq &
   \frac{-1}{C_1(\tilde M, \tilde B)}
   \sum_{v \in V}
     \frac{1-|\mu^v|^2}{1 + |\mu^v|^2}
 \sum_{e\in E_0(v)}
 c \, ( {\cal L}^e_0(t)  +  {\cal L}^e_1(t)   )
+
  4\,  \max_{e\in E} \nu^e \,\tilde B  \,
  \int_{t- T_0}^{t + T_0}
  {\cal L}_\Sigma(\tau)
   \, d\tau
\end{eqnarray*}
   
  We apply the observability inequalities \eqref{observabilityinequalityh1}  and obtain
 \begin{multline}\label{L0L1diff}
  {\cal L}_\Sigma(t+ T_0) -  {\cal L}_\Sigma(t - T_0)
 \\
 \leq 
 -
 \frac{c}{C_1(\tilde M, \tilde B)}
 \upsilon_0{\cal L}_\Sigma(t)
 + 
  4\,  \max_{e\in E} \nu^e \,\tilde B  \,
  \int_{t- T_0}^{t + T_0}
  {\cal L}_\Sigma(\tau)
   \, d\tau.
\end{multline}

   Now we need to control the integral on the right hand side of \eqref{L0L1diff}. To this end, we derive an
   estimate that shows
   that locally around $t$,
   the growth of
    ${\cal L}_\Sigma$ is
    limited.
   We have 
   \[
   \frac{d}{dt} {\cal L}_\Sigma(t)
   \leq
    4 \,  \nu^e
    \tilde B   \, 
    {\cal L}_\Sigma(t).
   \]
   Consider $\bar t \in (0, t)$.
   Then  for all $ s\in [ -\bar t, \, \bar t]$
   we have
  $
 {\cal L}_\Sigma(t + s)
\leq  
    \Delta(\nu \bar t) \,
 {\cal L}_\Sigma( t - \bar t ) 
 $
   with
  $ \Delta(\nu \bar t)
    =
  \exp\left( 8 \,\max_{e\in E} \nu^e
    \tilde B   \, \bar t \right).
    $
%
   In particular, we have
   \begin{equation}
   \label{upsilonungleichung}
    {\cal L}_\Sigma( t ) 
        \geq \exp\left(- 8 \,\max_{e\in E} \nu^e
    \tilde B   \,\bar t  \right)
 {\cal L}_\Sigma(t + \bar t )
  \\=  \frac{1}{\Delta(\nu \bar t)} \,  {\cal L}_\Sigma(t + \bar t ) .
\end{equation}       
Hence the increase of ${\cal L}_\Sigma$ is
limited.
Note that we did not show that
${\cal L}_\Sigma$ is decreasing.

Now we return to
the question of exponential decay of
 ${\cal L}_\Sigma$.
For this purpose in
order to be able to use the observability
inequality  we choose $\bar t =  T_0>0$. 

We apply \eqref{upsilonungleichung} in \eqref{L0L1diff} and obtain
\begin{multline*}
{\cal L}_\Sigma(t+ T_0) -  {\cal L}_\Sigma(t - T_0)
    \leq 
  -  \frac{c}{C_1( \tilde M, \, \tilde B)} 
  \upsilon_0\frac{1}{\Delta(\nu \, T_0)}
  \,
{\cal L}_\Sigma(t + T_0) 
\\
+
    8 \, T_0\,  \max_{e\in E} \nu^e  
    \, \tilde B \, \Delta(\nu \bar t) \, 
    {\cal L}_\Sigma(t  - T_0) .
\end{multline*}

Thus, we have
\begin{equation}
{\cal L}_\Sigma(t+ T_0) 
 \left[
 1
 +
 \frac{c  \upsilon_0 }{
 \Delta(\nu \,T_0) \, 
 C_1( {\tilde M},  {\tilde B})} 
 \right]
 \\
 \leq
 \left[
 1
  +
    8  \, T_0\,  \max_{e\in E} \nu^e  \,  \Delta(\nu \, T_0) \, 
 \right]
 {\cal L}_\Sigma  (t - T_0) .
\end{equation}
If
(\ref{observiassumption}) holds,
similarly as in  Lemma 2 from \cite{gugattucsnak}, this implies that ${\cal L}_\Sigma={\cal L}_0 +{\cal L}_1 $ decays exponentially fast.
Thus we have shown the exponential decay of the time derivatives.
 Thus we have proved Theorem   \ref{stability}.

\section{Numerical Experiments}
\label{Sec:Numerical}
In this section, we present  numerical experiments  to illustrate the theoretical results. 
For the discretization of the convective terms, we have used a finite difference upwind scheme in space and explicit Euler in time, while
the temporal discretization of the friction terms is implicit Euler. 
Therefore we can use the maximal time step allowed by the CFL condition so that discontinuities
in the solution are not smoothed out by numerical diffusion.

In order to verify the exponential convergence predicted by Theorem 6, 
we have plotted $\mathcal{L}_0$ (see  \eqref{Ledefinition}) over time for different values of the parameters $\mu^v$ (see definition of the systems ($\mathbf{S}$) and ($\mathbf{R}$)). For continuous solutions, we have also plotted $\mathcal{L}_1$ (see  \eqref{Ledefinitionobserv}).

\begin{figure}[h]
	\centering
	\includegraphics[width=0.7\linewidth]{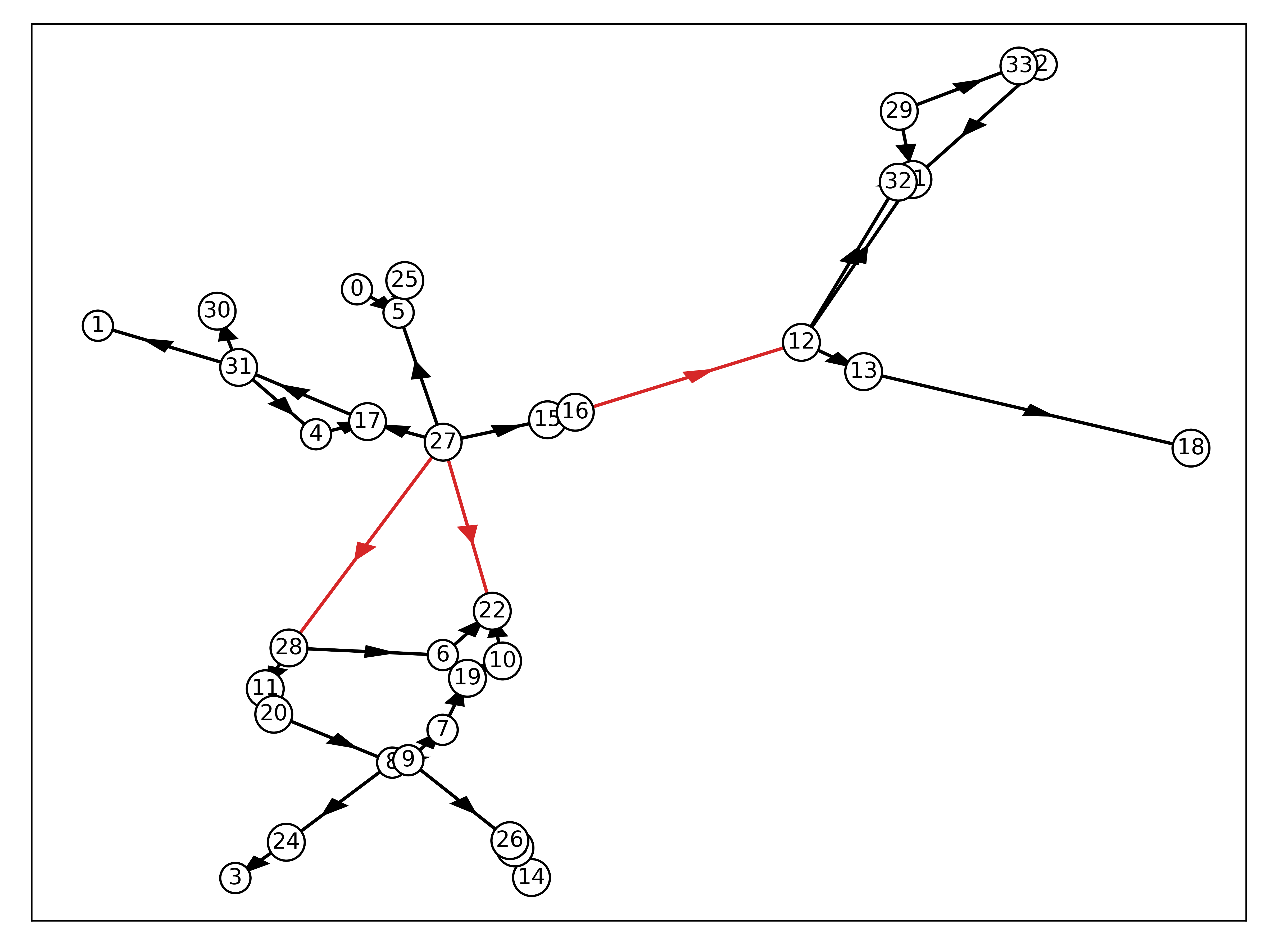}
	\caption{Sketch of the network for the numerical experiments. The arrows show the direction of the orientation of the pipes, which in general does not coincide with the direction of the flow. In red we have marked the edges on which the initial pressure  for the original system and the observer system differ.}
	\label{fig:networksketch}
\end{figure}

For the numerical experiments, we have used a modified version of the \linebreak GasLib-40 network (see [https://gaslib.zib.de],  \cite{data}; 
we have removed the compressors) that is shown in Figure~\ref{fig:networksketch}.
This network has 34 pipes with different lengths between \SI{3.068}{km} and \SI{86.690}{km} and diameters between \SI{0.4}{m} and \SI{1}{m}.
Additionally, we have used the parameters
$
	\theta = \SI{0.0137}{\frac{1}{m}}, \ c=\SI{340}{\frac{m}{s}}
	$
together with the pressure law $p(\rho)=c^2\rho$.
In all computations the initial velocity is  zero and both systems ($\mathbf{S}$) and ($\mathbf{R}$) have the same boundary value $u^e(t)$, which is piecewise linear in time
and can be computed 
for all $t>0$  
from the pressure
\begin{equation*}
    p=\begin{cases}\SI{59.5}{bar},&t=\SI{0}{s}\\ \SI{60.5}{bar},  & t=\SI{100}{s} \\ \SI{60}{bar}, &t\ge \SI{200}{s}
    \end{cases}
\end{equation*}
at all nodes and the (one-dimensional) mass flow
\begin{equation*}
    m=\begin{cases} \SI{41.788}{kg/s}, & \text{for nodes 0, 1 and 2 (`sources')}\\
    \SI{-4.323}{kg/s}, & \text{else}.
    \end{cases}
\end{equation*}

The numerical results confirm that
if  $|\mu^v|<1$ at all nodes, 
the difference $\delta$ between the system
($\mathbf{S}$) and the observer system ($\mathbf{R}$) decays to zero at least exponentially.

\subsection{Discontinuous initial data and
friction}
The initial data for the first experiment is a piecewise constant function, i.e., 
the initial pressure on the pipes connecting nodes 12 and 16, 27 and 28, 22 and 27 (see Figure \ref{fig:networksketch}) is $\SI{60}{bar}+h^e$ on the half of the pipes adjacent to the
node with smaller index and $\SI{60}{bar}$ on the other half of the pipe.
For the system~($\mathbf{S}$) we have used $h^e=\SI{2}{bar}$ for the edges connecting nodes 12 and 16, 27 and 28 and $h^e=\SI{1}{bar}$ for the edge between node 22 and 27,
while for the observer system~($\mathbf{R}$) we have used $h^e=\SI{1.5}{bar}$ and $h^e=\SI{0.75}{bar}$, respectively.
For all other pipes, the initial pressure is constant $\SI{60}{bar}$ for both systems.

\begin{figure}[h]
	\centering
	\includegraphics[width=0.5\linewidth]{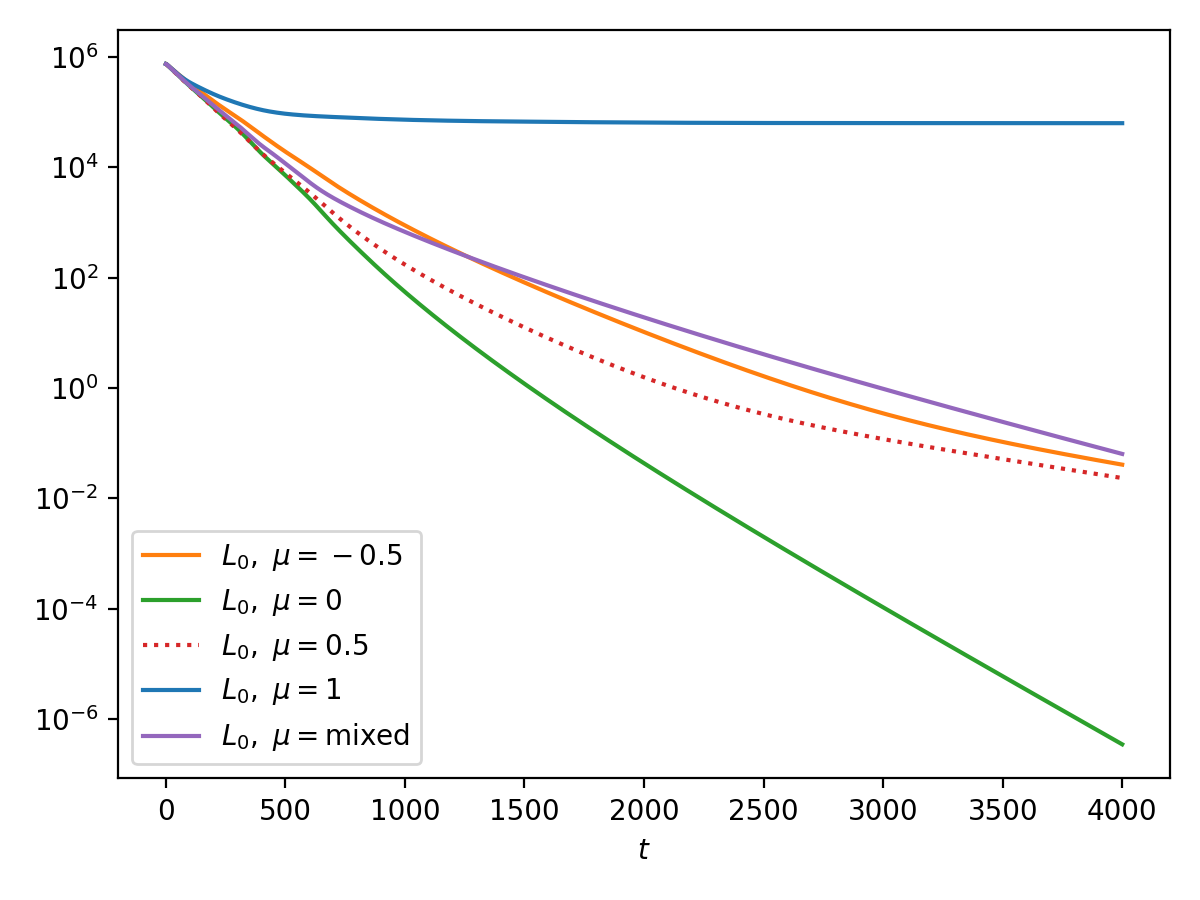}
	\caption{Discontinuous initial data and  friction: Temporal evolution of $\mathcal{L}_0$ for different values of $\mu^v$. In all but one cases we have set the same value of $\mu^v$ for all nodes, while `mixed' means that $\mu^v=0$ at
	all nodes with even index and additionally at the nodes 1, 5, 7, 15, 17, 29
	and $\mu^v=1$ for the remaining nodes.}
	\label{fig:L0stepfriction}
\end{figure}

\begin{figure}[h]
	\begin{subfigure}[c]{0.32\linewidth}
	\includegraphics[width=0.99\linewidth, trim=3.5cm 2cm 2.5cm 1cm, clip]{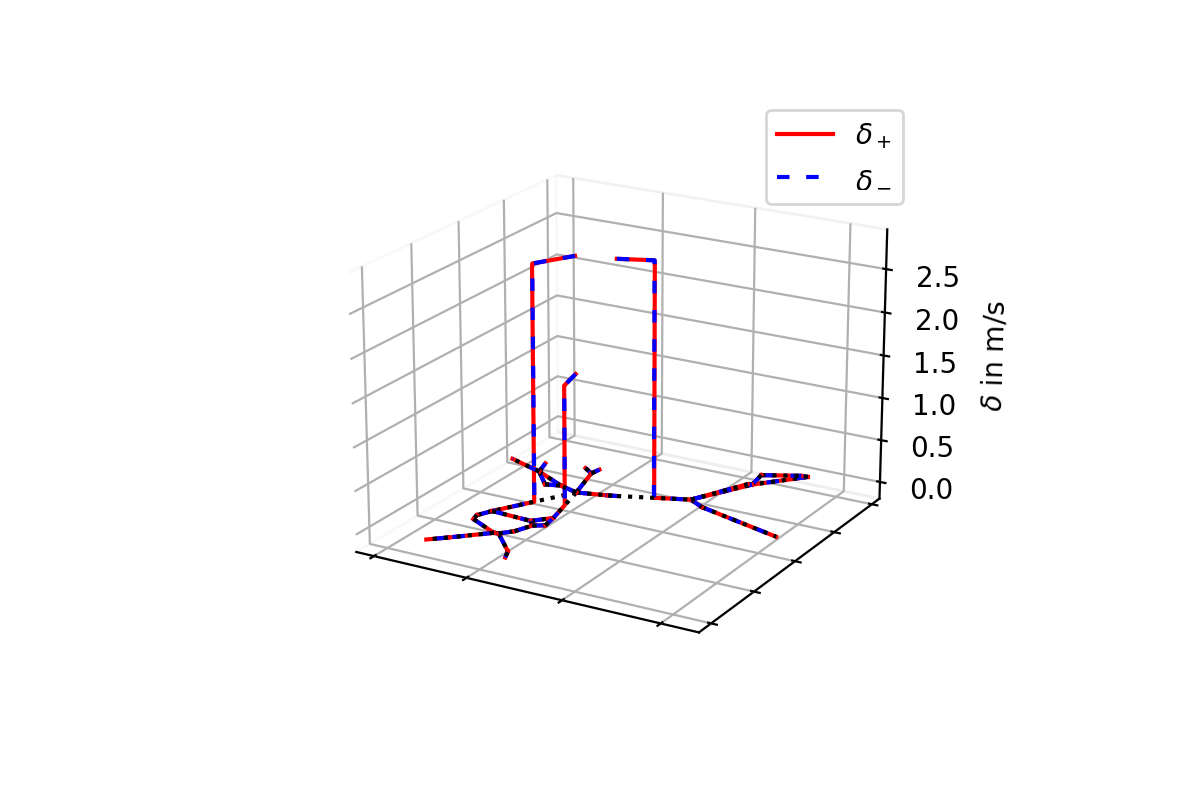}
	\end{subfigure}
	\begin{subfigure}[c]{0.32\linewidth}
	\includegraphics[width=0.99\linewidth, trim=3.5cm 2cm 2.5cm 1cm, clip]{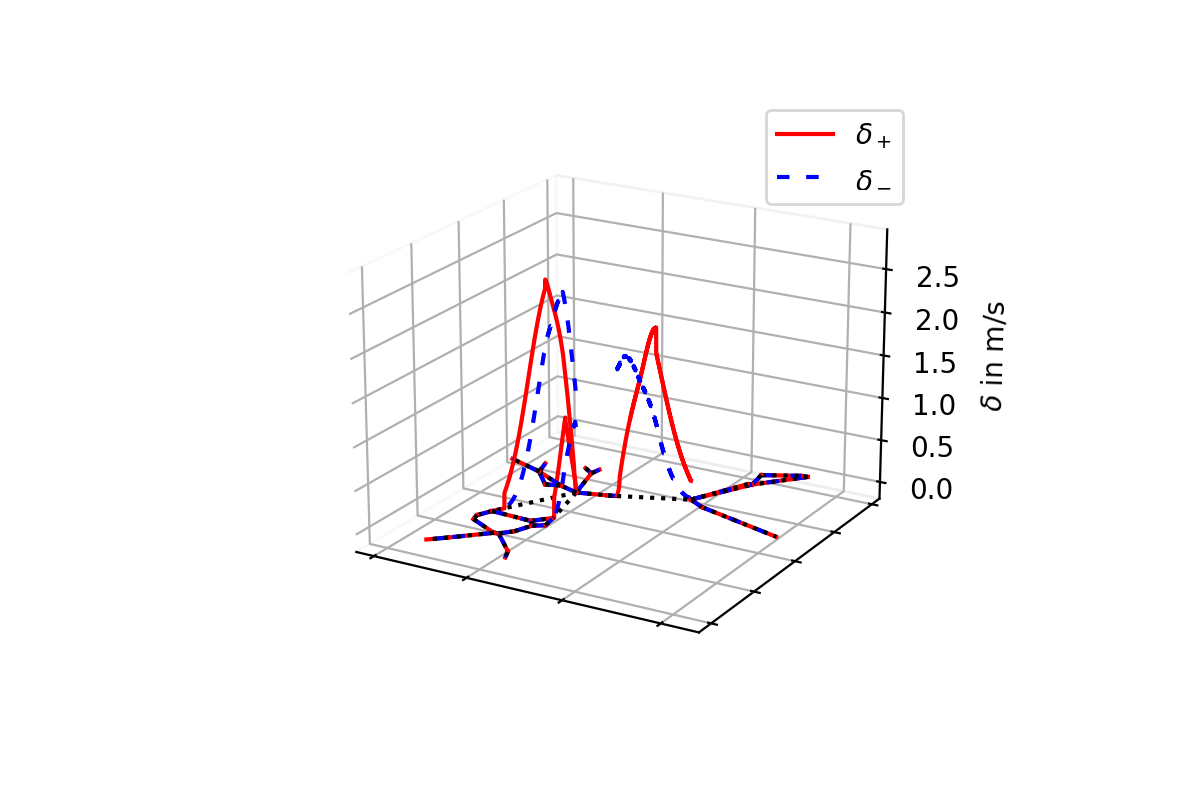}
	\end{subfigure}
	\begin{subfigure}[c]{0.32\linewidth}
	\includegraphics[width=0.99\linewidth, trim=3.5cm 2cm 2.5cm 1cm, clip]{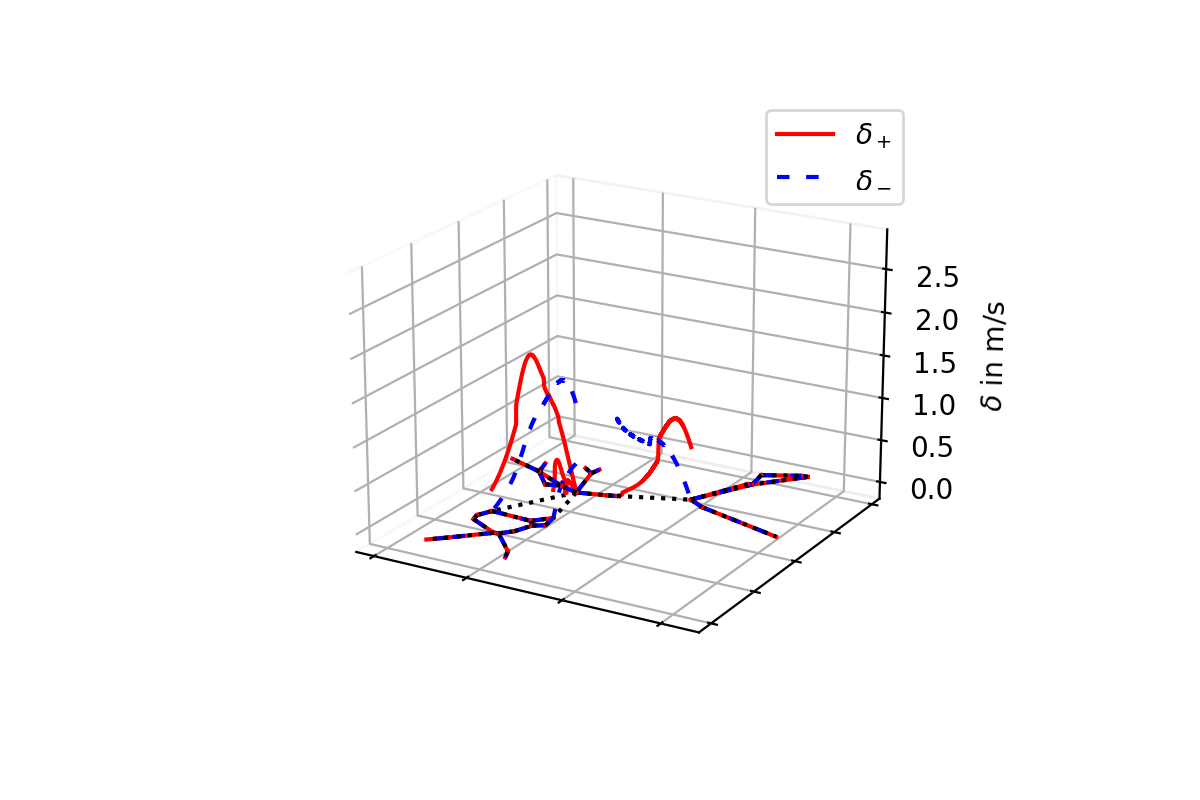}
	\end{subfigure}
	\caption{Discontinuous initial data and  friction: Snapshots of solutions at times $t=\SI{0}{s}$, $t\approx\SI{90}{s}$ and $t\approx\SI{180}{s}$ for $\mu^v=0$ at all nodes. On all pipes, $\delta_+$ is a red continuous line, $\delta_-$ is a blue dashed line and the network is shown by black dotted lines. }
	\label{fig:snapshots_step_friction}
\end{figure}

We plot the result for $\mathcal{L}_0$ in Figure \ref{fig:L0stepfriction}.
As predicted, $\mathcal{L}_0$ decays exponentially for all cases except for $\mu^v=1$ at all nodes, which is in accordance with the theoretical results.
As expected, we see the fastest convergence for $\mu^v=0$ at all nodes.
In addition, snapshots of the difference of the numerical solutions at times $t=\SI{0}{s}$, $t\approx\SI{90}{s}$ and $t\approx\SI{180}{s}$ for $\mu^v=0$ at all nodes are shown in Figure \ref{fig:snapshots_step_friction}.
The pictures show that the difference $\delta$ between the two systems decreases over time.

\subsection{Discontinuous initial data without friction}
\begin{figure}[h]
	\centering
	\includegraphics[width=0.5\linewidth]{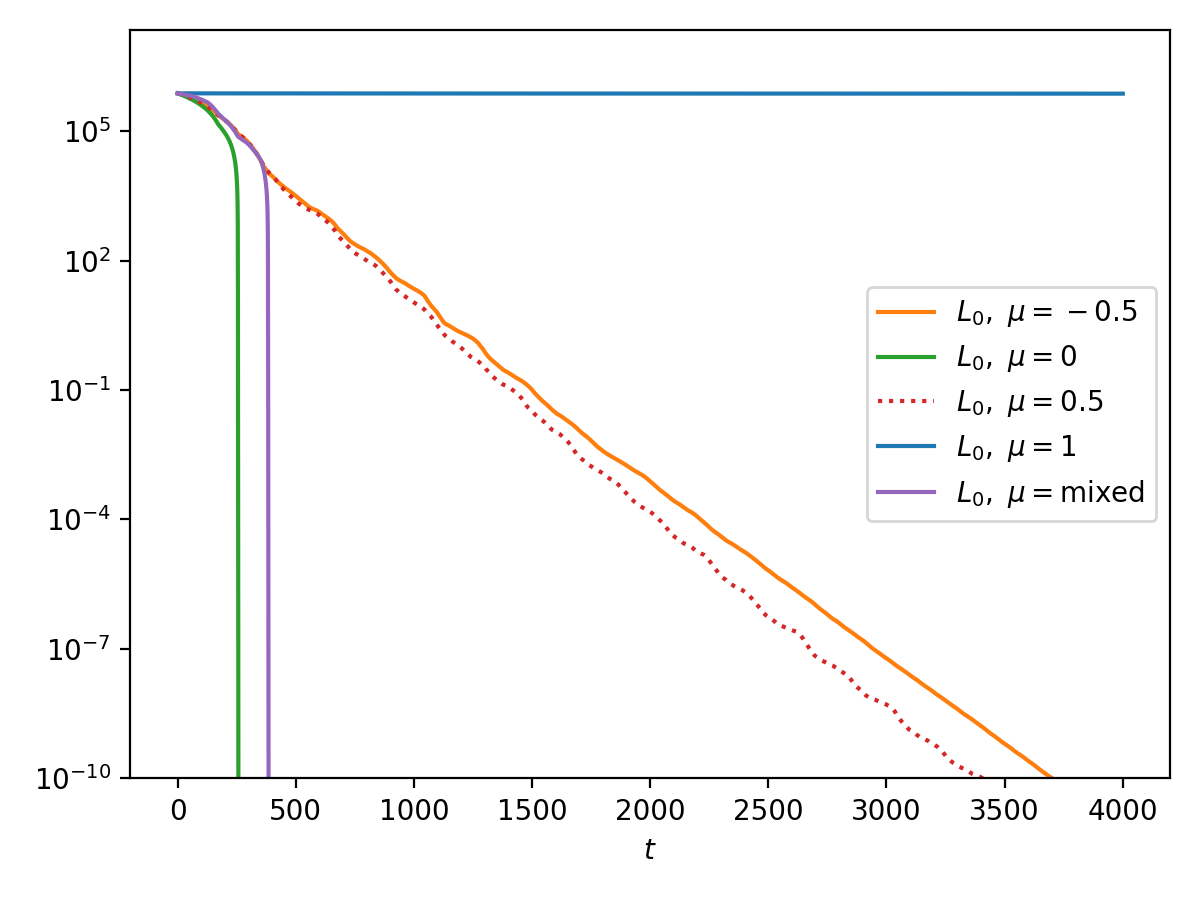}
	\caption{Discontinuous initial data without friction: Temporal evolution of $\mathcal{L}_0$ for different values of $\mu^v$. In all but one case we have set the same value of $\mu^v$ for all nodes, while `mixed' means that $\mu^v=0$ at
	all nodes with even index and additionally at the nodes 1, 5, 7, 15, 17, 29 and  $\mu^v=1$ for the remaining nodes.}
	\label{fig:L0stepnofriction}
\end{figure}

For this experiment we use the same initial data as in the previous section but set the  friction to zero (that is $\theta=0$).
We plot the result for $\mathcal{L}_0$ in Figure \ref{fig:L0stepnofriction}.
Again, $\mathcal{L}_0$ decays exponentially for $|\mu^v|<1$ at all nodes and for the `mixed' case. 
It can be seen that, in the case $\mu^v=0$ at all nodes, $\mathcal{L}_0$ vanishes after about $\SI{257.5}{s}$, which is approximately the time that a wave needs to travel through the longest pipe (that connects nodes 27 and 28 and has a length of about $\SI{86.7}{km}$).

For the case $\mu^v=0$ on half of the nodes and $\mu^v=1$ at the remaining nodes,  $\mathcal{L}_0$ vanishes after about $\SI{385}{s}$. 
So in the case of a source term that is zero,
the observer and the original system can be synchronized in finite-time.

\begin{figure}[h]
	\begin{subfigure}[c]{0.32\linewidth}
		\includegraphics[width=0.99\linewidth, trim=3.5cm 2cm 2.5cm 1cm, clip]{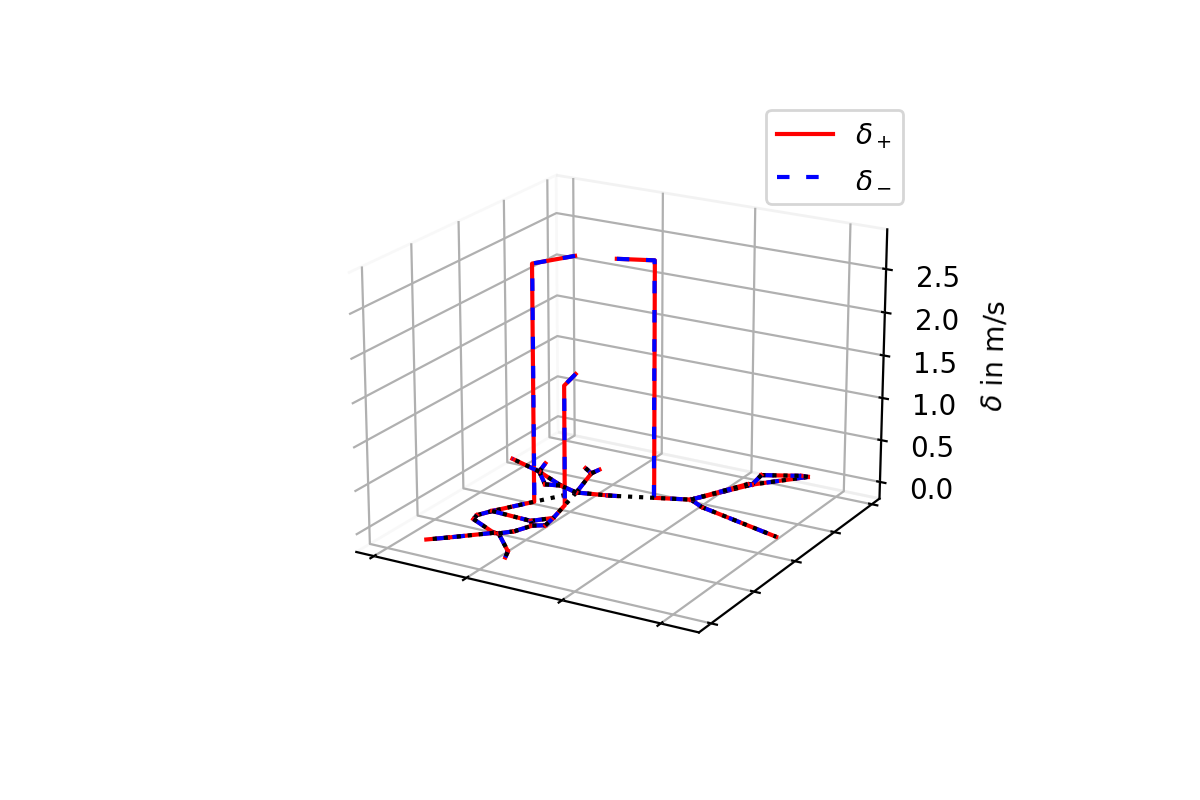}
	\end{subfigure}
	\begin{subfigure}[c]{0.32\linewidth}
		\includegraphics[width=0.99\linewidth, trim=3.5cm 2cm 2.5cm 1cm, clip]{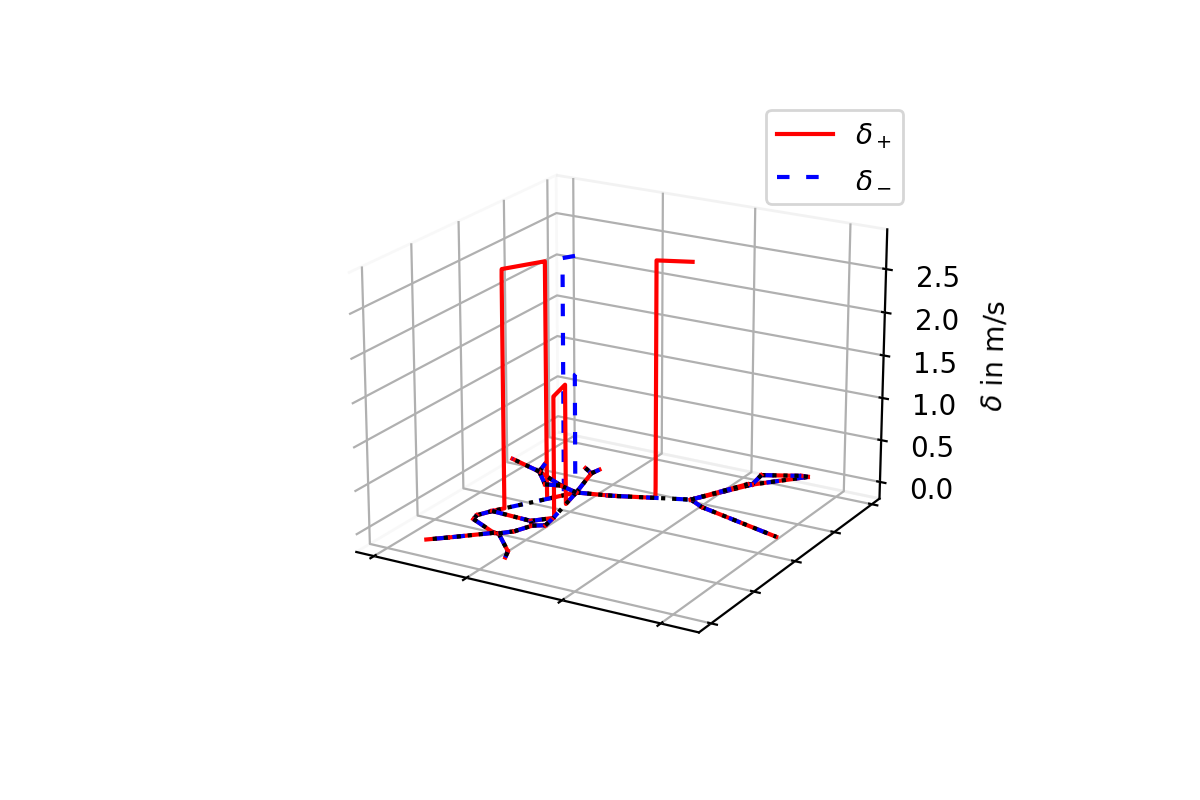}
	\end{subfigure}
	\begin{subfigure}[c]{0.32\linewidth}
		\includegraphics[width=0.99\linewidth, trim=3.5cm 2cm 2.5cm 1cm, clip]{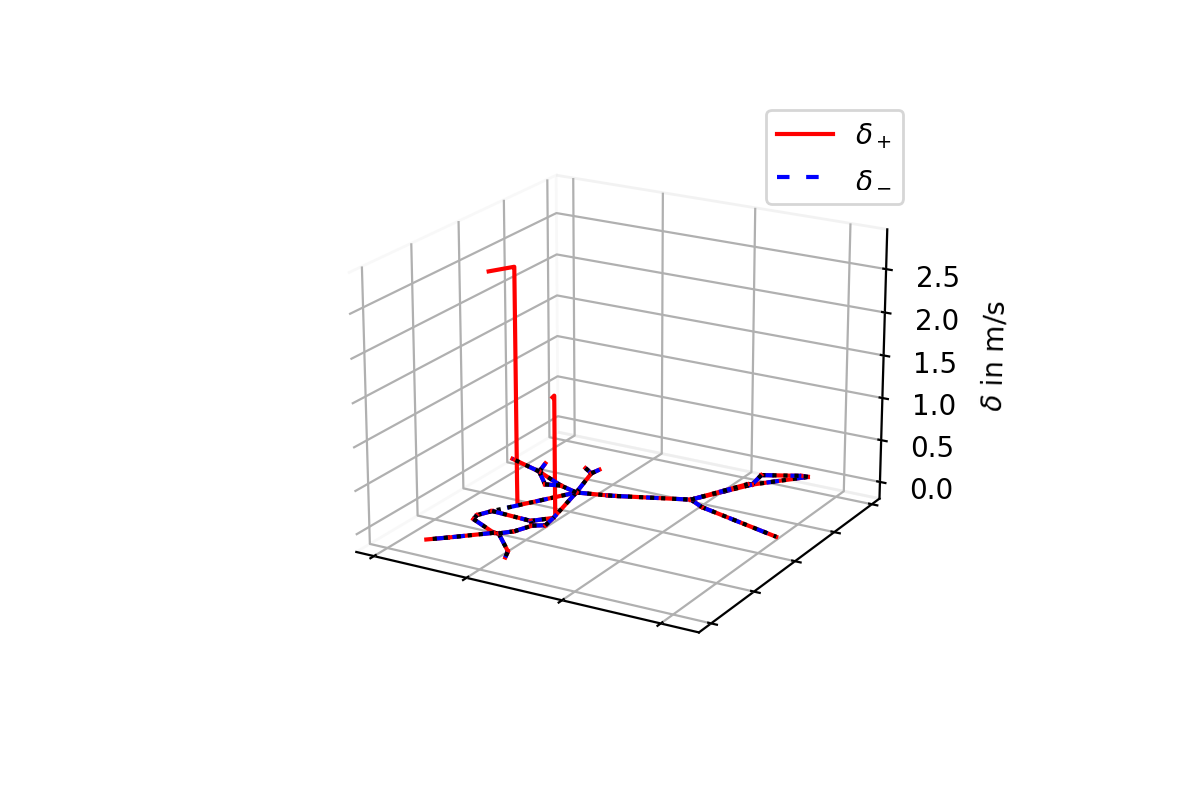}
	\end{subfigure}
	\caption{Discontinuous initial data without friction: Snapshots of solutions at times $t=\SI{0}{s}$, $t\approx\SI{90}{s}$ and $t\approx\SI{180}{s}$ for $\mu^v=0$ at all nodes. On all pipes, $\delta_+$ is a red continuous line, $\delta_-$ is a blue dashed line and the network is shown by black dotted lines. }
	\label{fig:snapshots_step_nofriction}
\end{figure}

In addition, snapshots of the difference of the numerical solutions at times $t=\SI{0}{s}$, $t\approx\SI{90}{s}$ and $t\approx\SI{180}{s}$ for $\mu^v=0$ at all nodes are shown in Figure~\ref{fig:snapshots_step_nofriction}.
The pictures show that the discontinuities in the solution remain,
since  the discretization of the convective terms produces no numerical diffusion.

\subsection{Continuous initial data and  friction}
In this experiment the initial pressure is a continuous function, i.e., 
the initial pressure on the pipes connecting nodes 12 and 16, 27 and 28, 22 and 27 (see Figure \ref{fig:networksketch}) is
${p(x)=\left(60+h^e \sin (f^e \frac{\pi x}{L})\right)\text{bar}}$.
For the system~($\mathbf{S}$) we have used $h^e=2$ and $f^e=2$ for the edges connecting nodes 12 and 16, 27 and 28 and $h^e=1$, $f^e=4$ for the edge between node 22 and 27,
while for the observer system~($\mathbf{R}$) we have used $h^e=1.5$, $f^e=2$ and $h^e=0.75$, $f^e=4$, respectively.
For all other pipes, the initial pressure is constant $\SI{60}{bar}$ for both systems.

\begin{figure}[h]
    \begin{subfigure}[c]{0.5\linewidth}
        \includegraphics[width=0.99\linewidth]{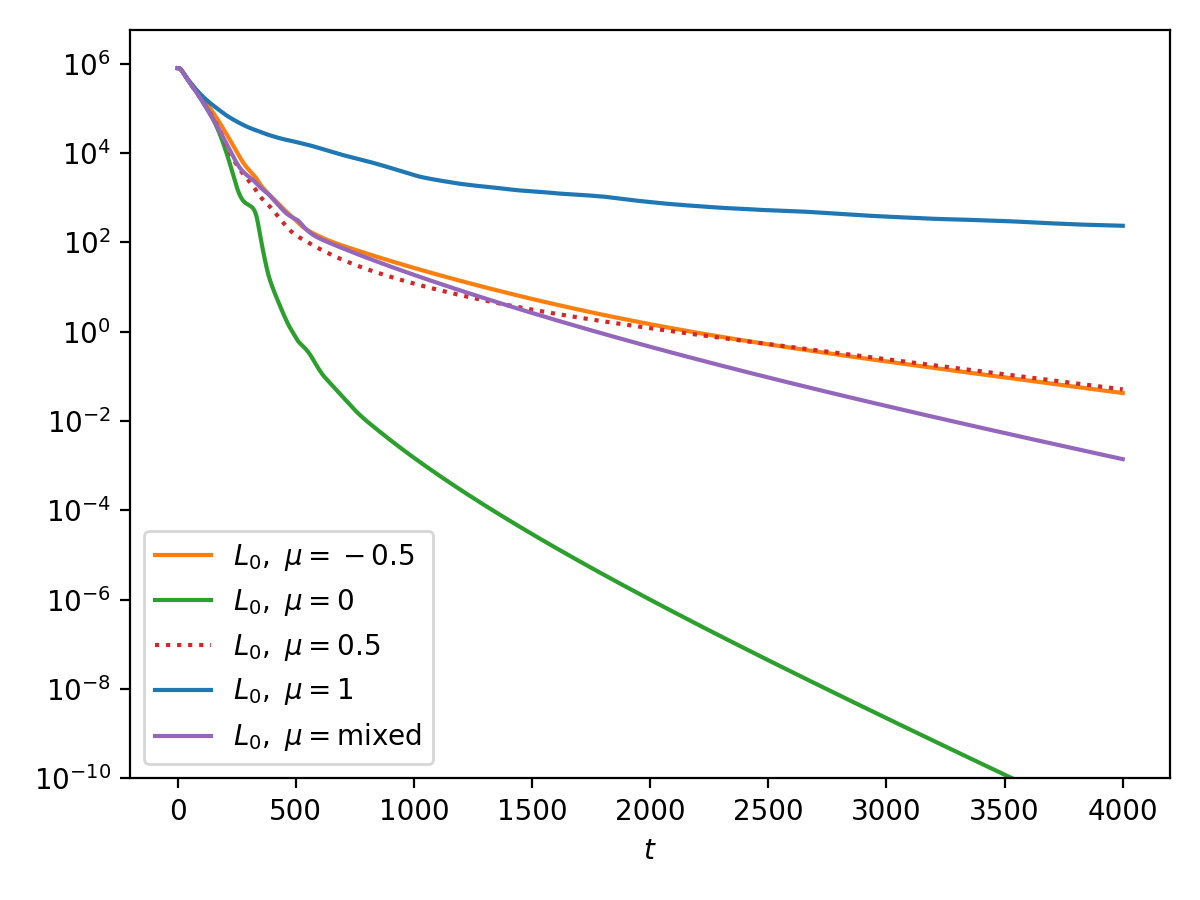}
    \end{subfigure}
    \begin{subfigure}[c]{0.5\linewidth}
        \includegraphics[width=0.99\linewidth]{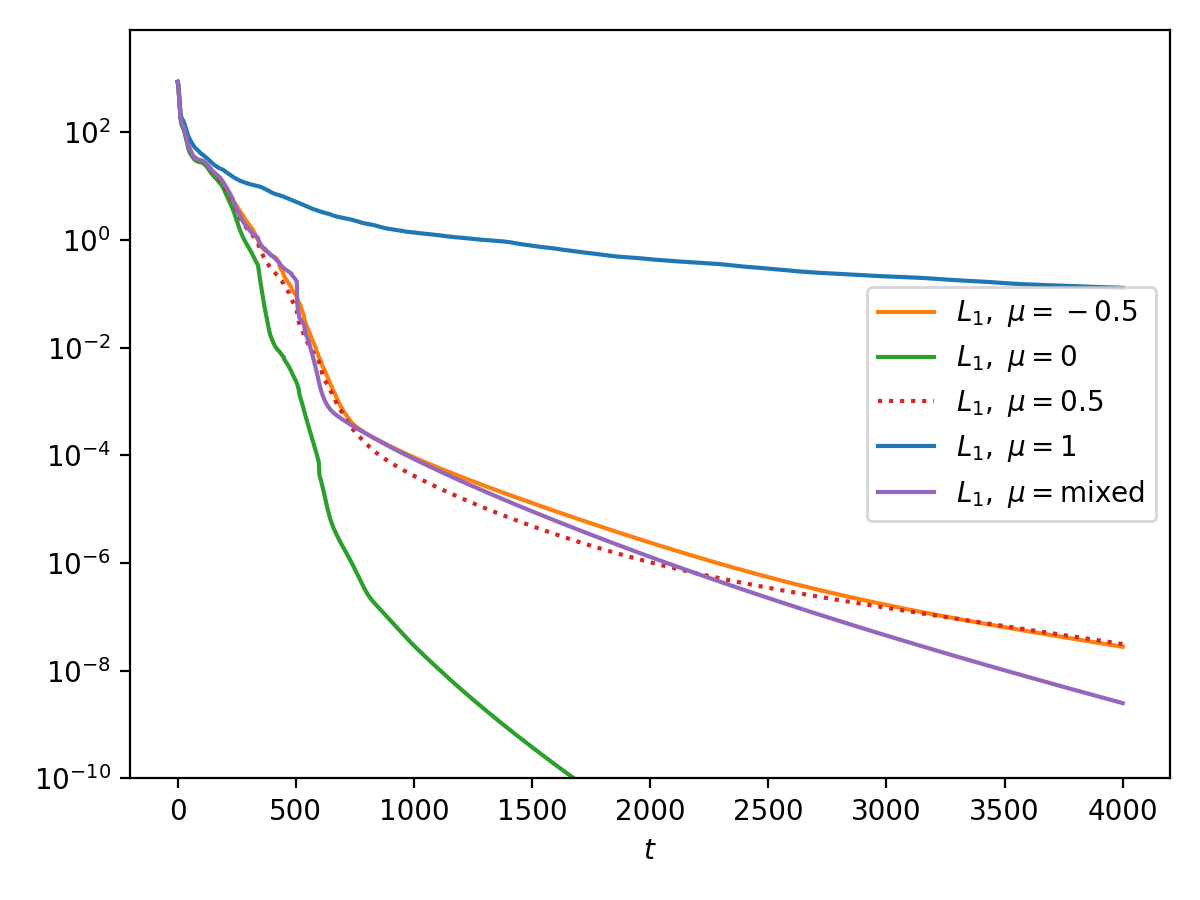}
    \end{subfigure}
	\caption{Continuous initial data and  friction: Temporal evolution of $\mathcal{L}_0$ (left) and $\mathcal{L}_1$ (right) for different values of $\mu^v$. In all but one case we have set the same value of $\mu^v$ for all nodes, while `mixed' means that $\mu^v=0$ at
	all nodes with even index and additionally at the nodes 1, 5, 7, 15, 17, 29 and  $\mu^v=1$ for the remaining nodes.}
	\label{fig:L0sinusfriction}
\end{figure}

We plot the results for $\mathcal{L}_0$  and $\mathcal{L}_1$ in Figure~\ref{fig:L0sinusfriction}.
As in the experiment with discontinuous initial data, $\mathcal{L}_0$ converges exponentially to zero for $|\mu^v|<1$ at all nodes and for the 'mixed' case.
Here we show additionally the value of $\mathcal{L}_1$, which also converges exponentially except in the case $\mu^v=1$ at all nodes.

\begin{figure}[h]
	\begin{subfigure}[c]{0.32\linewidth}
		\includegraphics[width=0.99\linewidth, trim=3.5cm 2cm 2.5cm 1cm, clip]{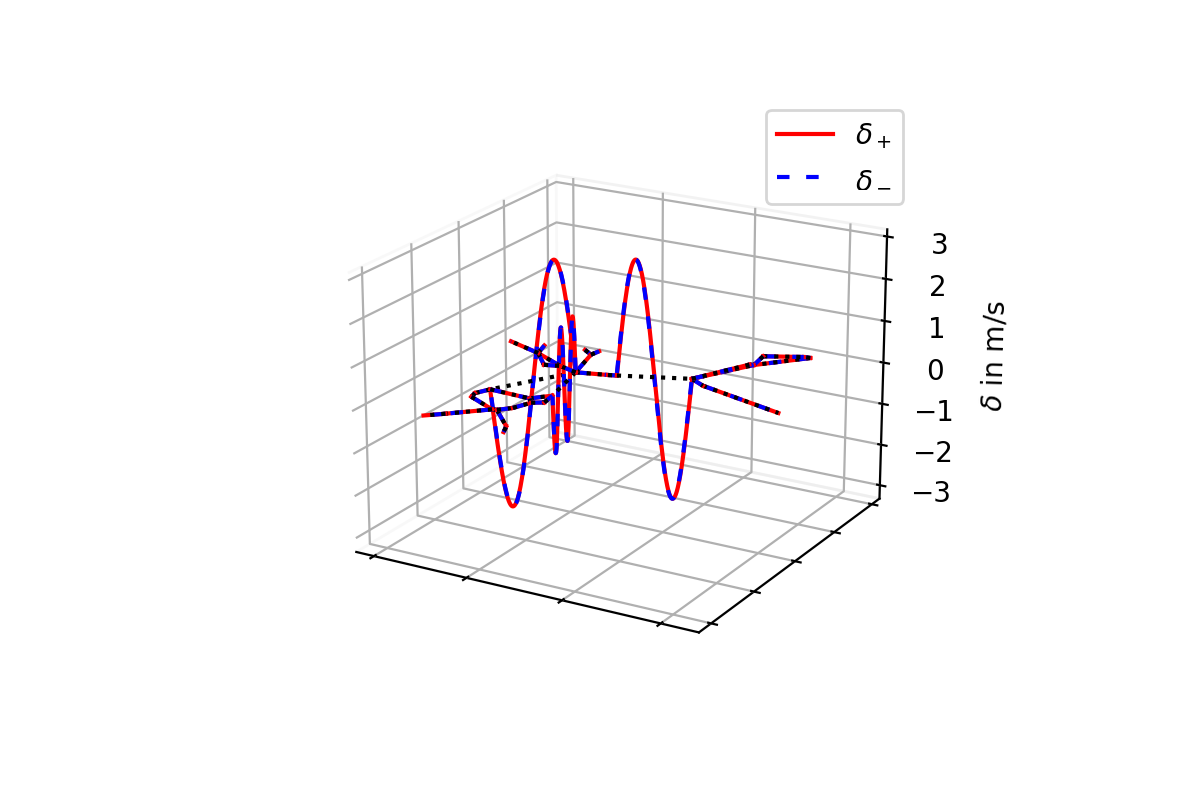}
	\end{subfigure}
	\begin{subfigure}[c]{0.32\linewidth}
		\includegraphics[width=0.99\linewidth, trim=3.5cm 2cm 2.5cm 1cm, clip]{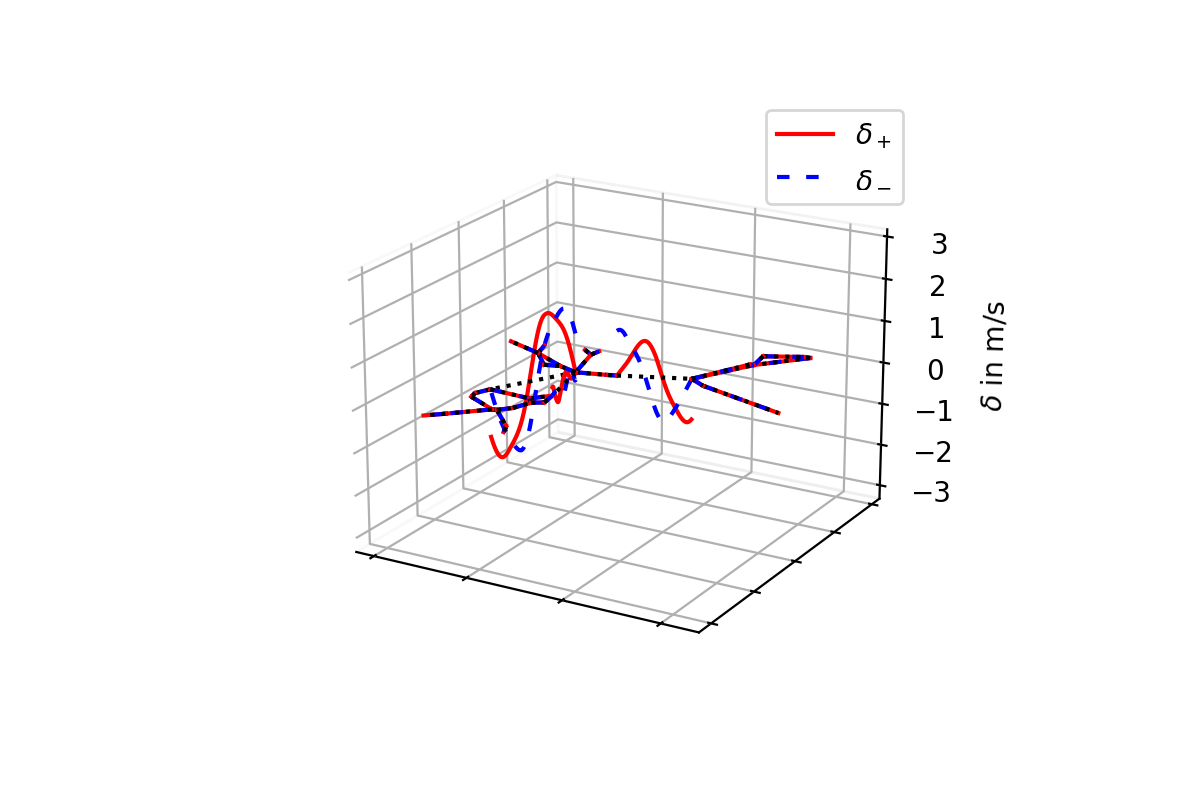}
	\end{subfigure}
	\begin{subfigure}[c]{0.32\linewidth}
		\includegraphics[width=0.99\linewidth, trim=3.5cm 2cm 2.5cm 1cm, clip]{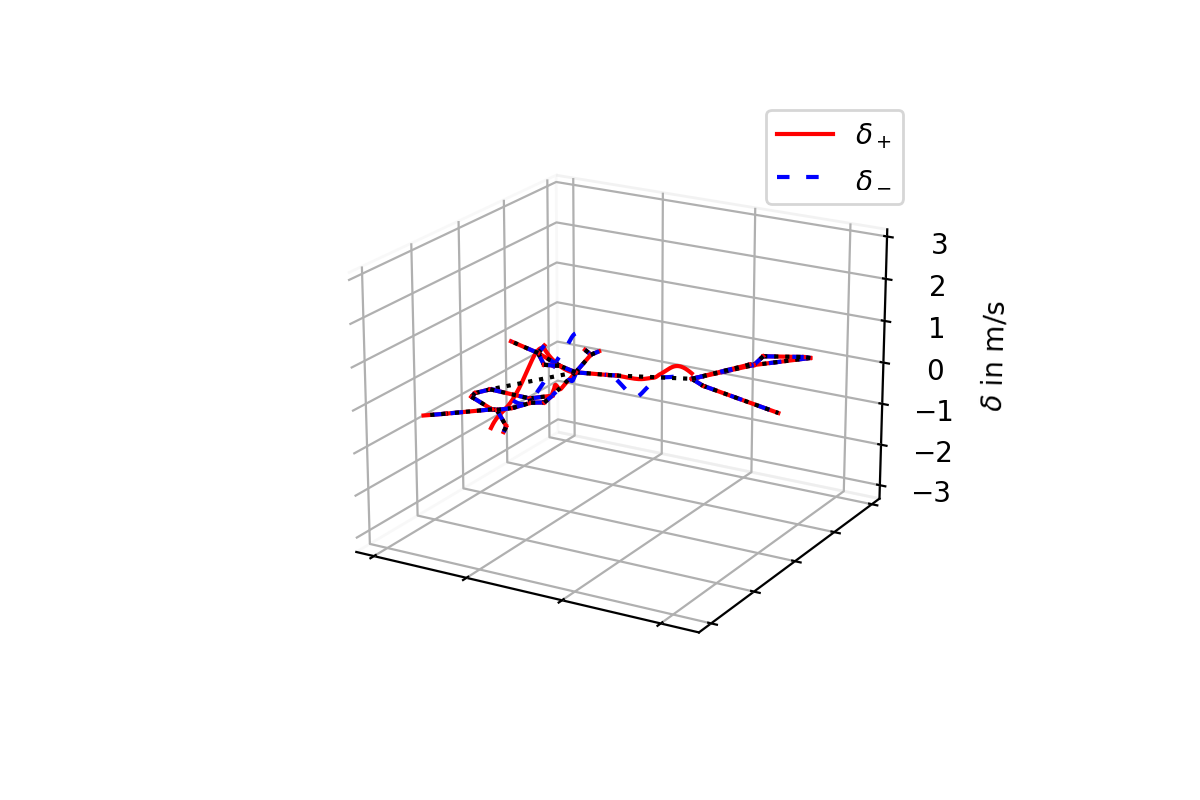}
	\end{subfigure}
	\caption{Continuous initial data and  friction: Snapshots of solutions at times $t=\SI{0}{s}$, $t\approx\SI{90}{s}$ and $t\approx\SI{180}{s}$ for $\mu^v=0$ at all nodes. On all pipes, $\delta_+$ is a red continuous line, $\delta_-$ is a blue dashed line and the network is shown by black dotted lines. }
	\label{fig:snapshots_sinus_friction}
\end{figure}

The snapshots of the difference of the numerical solutions at times $t=\SI{0}{s}$, $t\approx\SI{90}{s}$ and $t\approx\SI{180}{s}$ for $\mu^v=0$ at all nodes, displayed in Figure~\ref{fig:snapshots_sinus_friction},
show that the difference already reduces significantly in the first $\SI{180}{s}$.

\subsection{Continuous initial data without friction}
\begin{figure}[h]
    \begin{subfigure}[c]{0.5\linewidth}
        \includegraphics[width=0.99\linewidth]{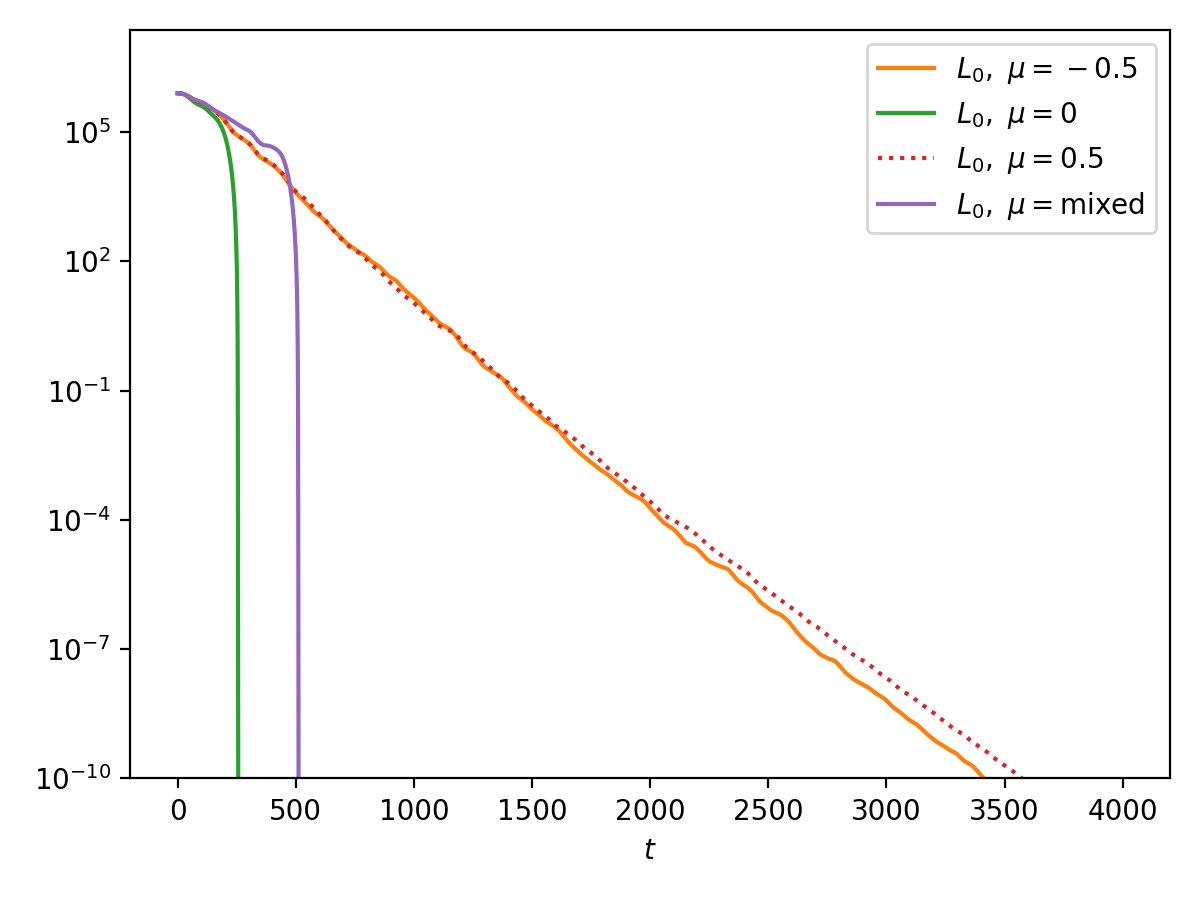}
    \end{subfigure}
    \begin{subfigure}[c]{0.5\linewidth}
        \includegraphics[width=0.99\linewidth]{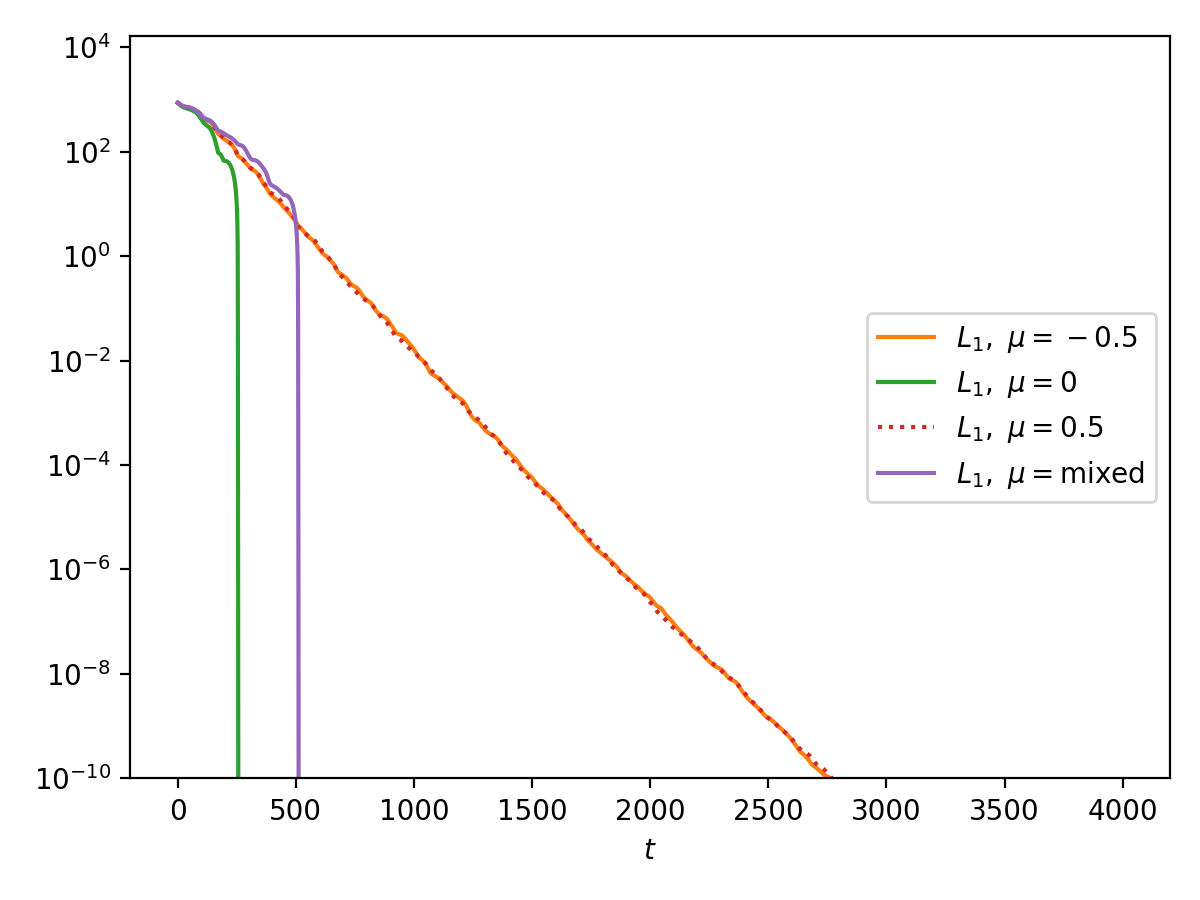}
    \end{subfigure}
	\caption{Continuous initial data without friction: Temporal evolution of $\mathcal{L}_0$ (left) and $\mathcal{L}_1$ (right) for different values of $\mu^v$. In all but one case we have set the same value of $\mu^v$ for all nodes, while `mixed' means that $\mu^v=0$ at
	all nodes with even index and additionally at the nodes 1, 5, 7, 15, 17, 29 and  $\mu^v=1$ for the remaining nodes.}
	\label{fig:L0sinusnofriction}
\end{figure}

Now we consider the case without friction, while we use the same initial data as in the previous experiment.
We plot the results for $\mathcal{L}_0$ and 
$\mathcal{L}_1$ in Figure \ref{fig:L0sinusnofriction}.
Similar to the experiment for discontinuous initial data without friction,
$\mathcal{L}_0$ and $\mathcal{L}_1$ vanish after about $\SI{257.5}{s}$ for 
$\mu^v=0$ at all nodes and after about $\SI{512.5}{s}$ for the `mixed' case.

\begin{figure}[h]
	\begin{subfigure}[c]{0.32\linewidth}
		\includegraphics[width=0.99\linewidth, trim=3.5cm 2cm 2.5cm 1cm, clip]{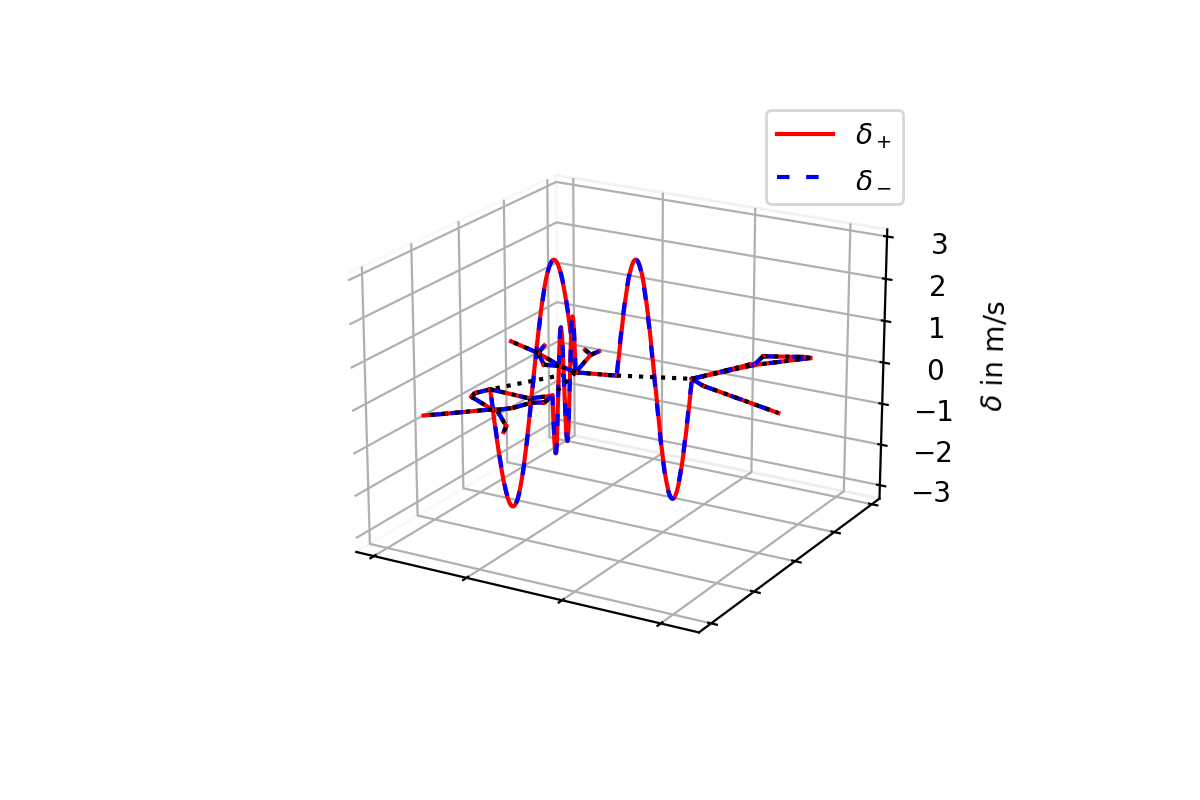}
	\end{subfigure}
	\begin{subfigure}[c]{0.32\linewidth}
		\includegraphics[width=0.99\linewidth, trim=3.5cm 2cm 2.5cm 1cm, clip]{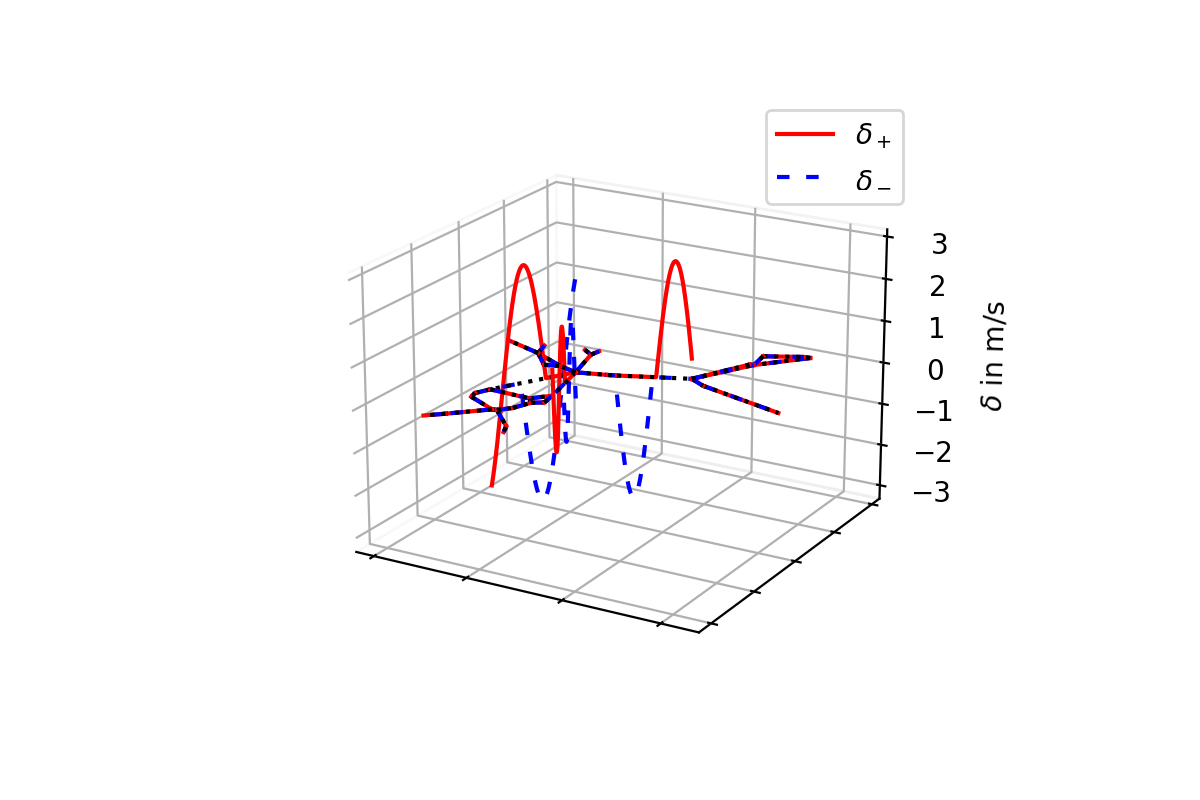}
	\end{subfigure}
	\begin{subfigure}[c]{0.32\linewidth}
		\includegraphics[width=0.99\linewidth, trim=3.5cm 2cm 2.5cm 1cm, clip]{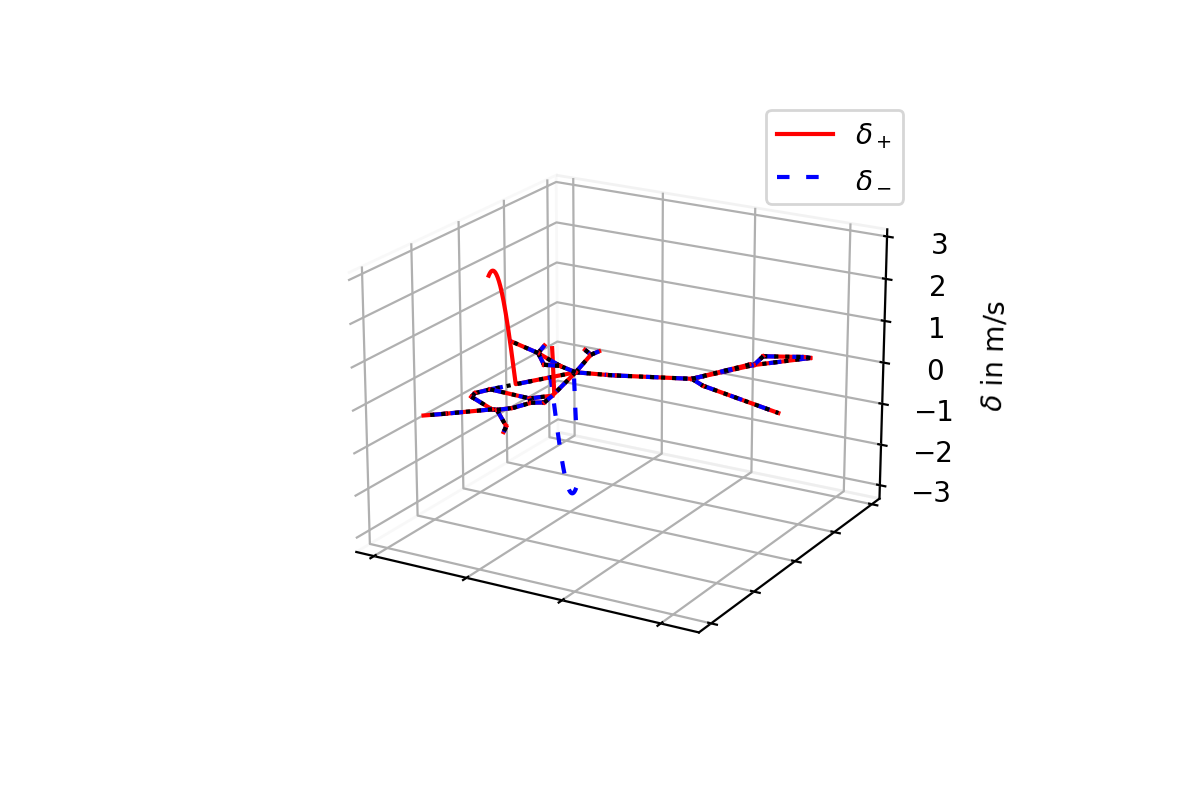}
	\end{subfigure}
	\caption{Continuous initial data without friction: Snapshots of solutions at times $t=\SI{0}{s}$, $t\approx\SI{90}{s}$ and $t\approx\SI{180}{s}$ for $\mu^v=0$ at all nodes. On all pipes, $\delta_+$ is a red continuous line, $\delta_-$ is a blue dashed line and the network is shown by black dotted lines. }
	\label{fig:snapshots_sinus_nofriction}
\end{figure}

The snapshots of the difference of the numerical solutions at times $t=\SI{0}{s}$, $t\approx\SI{90}{s}$ and $t\approx\SI{180}{s}$ for $\mu^v=0$ at all nodes in Figure~\ref{fig:snapshots_sinus_nofriction}
show that the difference is zero on large parts of the network after $\SI{180}{s}$.

\section{Conclusion}
In this paper, we have analyzed the performance of an observer system for
the gas flow through a pipeline network 
that is governed by a semilinear model.
As input data for the  observer system,
measurement data that is  obtained at certain 
points in space in  the network is used.
We have shown that under suitable
regularity conditions
for the solution
the observation error decays exponentially.
The theoretical findings are illustrated
by numerical experiments.

\bibliographystyle{siam}





\end{document}